\documentclass[11pt,reqno]{amsart}
\usepackage[utf8]{inputenc}
\usepackage{amssymb, mathtools}
\usepackage{tikz-cd}
\usepackage{dsfont}        %
\usepackage{newtxtext}
\usepackage{newtxmath}      %
\usepackage{mathrsfs}      %
\usepackage{stmaryrd}      %
\DeclareFontFamily{U}{rsfs}{\skewchar\font127 }
\DeclareFontShape{U}{rsfs}{m}{n}{<-6.5> rsfs5 <6.5-8> rsfs7 <8-> rsfs10}{}
\DeclareFontFamily{U}{stmry}{}
\DeclareFontShape{U}{stmry}{m}{n}{<-5.5> stmary5 <5.5-6.5> stmary6 <6.5-7.5> stmary7
  <7.5-8.5> stmary8 <8.5-9.5> stmary9 <9.5-> stmary10}{}

\usepackage{xcolor}
\definecolor{themegreen}{RGB}{0,92,54}   %
\definecolor{themeorange}{RGB}{191,87,0} %
\usepackage[final, colorlinks=true, linkcolor=themegreen, citecolor=themegreen, urlcolor=themegreen,pagebackref]{hyperref}
\usepackage[nameinlink]{cleveref}

\usepackage{tikz}

\usepackage{enumitem}
\usepackage[margin=28mm,a4paper]{geometry}
\usepackage[indent]{parskip}

\theoremstyle{plain}
\newtheorem{theorem}{Theorem}[section]
\newtheorem{introtheorem}{Theorem}

\newtheorem*{introcorollary}{Corollary}
\newtheorem{proposition}[theorem]{Proposition}
\newtheorem{lemma}[theorem]{Lemma}
\newtheorem{corollary}[theorem]{Corollary}

\theoremstyle{definition}
\newtheorem{definition}[theorem]{Definition}
\newtheorem{example}[theorem]{Example}
\newtheorem{remark}[theorem]{Remark}

\crefname{introtheorem}{Theorem}{Theorems}
\Crefname{introtheorem}{Theorem}{Theorems}
\crefname{question}{Question}{Questions}
\Crefname{question}{Question}{Questions}
\crefname{conjecture}{Conjecture}{Conjectures}
\Crefname{conjecture}{Conjecture}{Conjectures}

\makeatletter
\newif\ifshowlabelkeys
\@ifclasswith{amsart}{draft}{\showlabelkeystrue}{}
\@for\statementenv:=theorem,introtheorem,proposition,lemma,corollary,conjecture,definition,example,question,remark\do{%
  \@namedef{showstatementlabel@\statementenv}{}%
}
\ifshowlabelkeys
  \AtBeginDocument{%
    \let\statementoriginallabel\label
    \RenewDocumentCommand{\label}{om}{%
      \IfNoValueTF{#1}{\statementoriginallabel{#2}}{\statementoriginallabel[#1]{#2}}%
      \ifcsname showstatementlabel@\@currenvir\endcsname
        \marginpar{\raggedright\normalfont\tiny\ttfamily\color{black!60}%
          \hyphenchar\font=`\-\relax \hyphenpenalty=10000 \exhyphenpenalty=0
          \detokenize{#2}\par}%
      \fi
    }%
  }
\fi
\makeatother

\newcommand{\bk}{k}                           %
\newcommand{\PP}{\mathds{P}}

\newcommand{\ZZ}{\mathds{Z}}
\newcommand{\Fp}{\mathds{F}_p}
\newcommand{\Oh}{\mathcal{O}}

\newcommand{\Frob}{F}                           %
\DeclareMathOperator{\Spec}{Spec}
\DeclareMathOperator{\Hom}{Hom}

\DeclareMathOperator{\NS}{NS}

\DeclareMathOperator{\HH}{H}

\DeclareMathOperator{\Pic}{Pic}

\DeclareMathOperator{\rank}{rank}

\newcommand{\GGm}{\mathds{G}_m}

\newcommand{\QQ}{\mathds{Q}}

\DeclareMathOperator{\disc}{disc}
\numberwithin{equation}{section}

\newcommand{\A}{\mathds{A}}
\newcommand{\FF}{\mathds{F}}

\newcommand{\GGa}{\mathds{G}_a}
\newcommand{\cE}{\mathcal{E}}

\DeclareMathOperator{\Br}{Br}

\DeclareFontFamily{U}{wncy}{}
\DeclareFontShape{U}{wncy}{m}{n}{<->wncyr10}{}
\DeclareSymbolFont{cyrletters}{U}{wncy}{m}{n}
\DeclareMathSymbol{\Sha}{\mathord}{cyrletters}{"58}
\DeclareMathOperator{\Gal}{Gal}
\DeclareMathOperator{\Aut}{Aut}

\DeclareMathOperator{\cha}{char}

\newcommand{\sep}{\mathrm{sep}}

\newcommand{\TT}{\mathds{T}}                   %
\newcommand{\cG}{\mathcal{G}}                  %
\newcommand{\cT}{\mathcal{T}}                  %
\newcommand{\bH}{\mathbf{H}}                   %
\newcommand{\mup}{\mu_p}
\newcommand{\kab}{\overline{\bk(a)}}           %
\newcommand{\pr}{\mathrm{pr}}                  %
\DeclareMathOperator{\ord}{ord}                %
\newcommand{\KI}{\mathrm{I}}
\newcommand{\KII}{\mathrm{II}}
\newcommand{\KIII}{\mathrm{III}}
\newcommand{\KIV}{\mathrm{IV}}

\hypersetup{pdftitle={Rational multisections and unirationality of supersingular K3 surfaces}, pdfauthor={Haitao Zou}}
\title[Rational multisections and unirationality]{Rational multisections and unirationality of supersingular K3 surfaces}
\author{Haitao Zou}
\address{Fakult\"at f\"ur Mathematik, Universit\"at Bielefeld, 33615 Bielefeld, Germany}
\email{hzou@math.uni-bielefeld.de}
\date{December 1, 2026}
\thanks{The author is supported by the Deutsche Forschungsgemeinschaft (DFG, German Research Foundation), Project-ID 491392403, TRR 358.}
\AddToHook{env/proposition/begin}{\crefalias{theorem}{proposition}}
\AddToHook{env/lemma/begin}{\crefalias{theorem}{lemma}}
\AddToHook{env/corollary/begin}{\crefalias{theorem}{corollary}}
\AddToHook{env/conjecture/begin}{\crefalias{theorem}{conjecture}}
\AddToHook{env/definition/begin}{\crefalias{theorem}{definition}}
\AddToHook{env/example/begin}{\crefalias{theorem}{example}}
\AddToHook{env/remark/begin}{\crefalias{theorem}{remark}}
\AddToHook{env/question/begin}{\crefalias{theorem}{question}}
\begin{document}
\begin{abstract}
We prove that every supersingular $K3$ surface of Artin invariant at most nine is unirational, confirming the question of Artin, Rudakov, Shafarevich and Shioda outside the generic case of Artin invariant ten.
In odd characteristic, we reduce the question to Artin invariant nine and sweep out the $K3$ surface by Abel--Jacobi images of rational multisections of its Artin--Tate twists which vary along the supersingular twistor line.
\end{abstract}
\maketitle
\section{Introduction}
Let $k$ be an algebraically closed field of characteristic $p>0$.
A $K3$ surface $X$ over $k$ is \emph{(Shioda)-supersingular} if its Picard number is $\rho(X)=22$.
The discriminant of its N\'eron--Severi lattice $\NS(X)$ is then $-p^{2\sigma}$ for an integer $1\leq\sigma\leq10$, the \emph{Artin invariant} $\sigma(X)$ of $X$.
The locus of supersingular $K3$ surfaces of Artin invariant at most $\sigma$ has dimension $\sigma-1$ in odd characteristic \cite{OgusTorelli} and in characteristic two \cite{RSchar2}.
In particular, the general supersingular $K3$ surface has Artin invariant ten.

Shioda \cite{ShiodaUnirational} proved that the Fermat quartic surface is unirational as a $K3$ surface when $p \equiv 3 \pmod 4$. Namely, such a $K3$ surface admits a dominant rational map from $\PP^2$ which is purely inseparable. He also observed that unirational $K3$ surfaces are supersingular in general.
A conjecture of Artin, Rudakov, Shafarevich and Shioda predicts the converse: every supersingular $K3$ surface is unirational.
This is known in characteristic two \cite{RSchar2}; in odd characteristic for $\sigma\leq2$, where the surfaces are Kummer surfaces \cite{ShiodaKummer,Ogus}; in characteristic three for $\sigma\leq6$ \cite{RSchar2}; and in characteristic five for $\sigma\leq3$ \cite{PhoShimada}.

A proof for $p\geq5$ was proposed in \cite{Lie} by establishing a Shioda--Inose type sandwich theorem in positive characteristic, but it is incomplete \cite{erratum}: it relies on the existence of purely inseparable multisections on certain elliptic fibrations, which was disproved by Bragg and Lieblich \cite{BLpoints}.
Our main result establishes the conjecture outside the open stratum of Artin invariant ten with a completely different argument.

\begin{introtheorem}\label{thm:main}
    Every supersingular $K3$ surface of Artin invariant $\sigma\leq9$ is unirational.
\end{introtheorem}
The proof of \Cref{thm:main} establishes the unirationality of surfaces of Artin invariant nine directly. The smaller Artin invariants are deduced from it by \Cref{prop:AT twistor family}. 
In particular, the proof does not rely on Shioda's unirationality theorem on Kummer surfaces.

The unirationality of Enriques surfaces is closed related to that of $K3$ surfaces. 
In characteristic two, there is a complete classification by Blass and Crew; see \cite[Theorem 2]{Blass}.
In odd characteristic, Jang \cite[Corollary~2.4]{Jang15} shows that if the $K3$ cover of an Enriques surface is supersingular, its Artin invariant is at most $5$. By Serre's criterion \cite[Lemma~3.1]{ShiodaKummer}, \Cref{thm:main} has the following consequence.
\begin{introcorollary}
Suppose that $p>3$. An Enriques surface over $k$ is unirational if and only if its $K3$ cover is supersingular.
\end{introcorollary}

We also obtain partial results for Artin invariant ten, while the general case remains open in odd characteristic.

\begin{introtheorem}\label{thm:B}
 Let $X$ be a supersingular $K3$ surface of Artin invariant $10$.
 Suppose that $X$ admits an elliptic fibration with supersingular geometric generic fiber.
 Then $X$ is unirational.
\end{introtheorem}
Here, we say that a non-Jacobian elliptic fibration on a $K3$ surface has \emph{supersingular geometric generic fiber} if the Jacobian of its geometric generic fiber is a supersingular elliptic curve.
In this case we construct a purely inseparable dominant rational map from its Jacobian to the surface.
This allows us to deduce unirationality from the case of lower Artin invariant (\Cref{lem:purely inseparable sandwich}). Bragg and Lieblich's purely inseparable multisection result \cite{BLpoints} gives another proof in this setting, as explained in \Cref{rem:BL22-multisection}.

For $p\geq5$, the assumption of \Cref{thm:B} can be satisfied only when $p\equiv2\pmod3$, by \Cref{pencils:large}. For the case $p \equiv 1 \pmod 3$, there is an example of supersingular $K3$ surface of Artin invariant $10$, which is dominant by a unirational Fermat surface (\Cref{ex:ordinary Fermat}).
\subsection*{Strategy of the proof}
Let $p$ be odd, and let $J$ be a supersingular $K3$ surface with $\sigma(J)\leq9$.
It has a genus-one fibration $J\to\PP^1$ with a section (\Cref{prelim:jacobian}), which we may assume to be elliptic, since a quasi-elliptic surface is unirational \cite[\S2, Proposition~2]{RS}. Let $E/k(t)$ be its generic fiber.
The Brauer group $\Br(J)$, which is the Tate--Shafarevich group of $E$, is nonzero and killed by $p$.
The torsor of a nonzero class $\alpha\in\Br(J)$ has a relatively minimal regular model $X(\alpha)\to\PP^1$, which is a supersingular $K3$ surface of Artin invariant $\sigma(J)+1$; we call it an \emph{Artin--Tate twist} of $J$ (\Cref{prelim:artin-tate}).
Since $p\alpha=0$, multiplication by $p$ induces a dominant rational map $X(\alpha)\dashrightarrow J$.
Hence unirationality in Artin invariant $\sigma+1$ implies unirationality in Artin invariant $\sigma$ (\Cref{prop:AT twistor family}), and it suffices to treat $\sigma(J)=9$.

\begin{figure}[t]
  \centering
  \includegraphics[width=\textwidth]{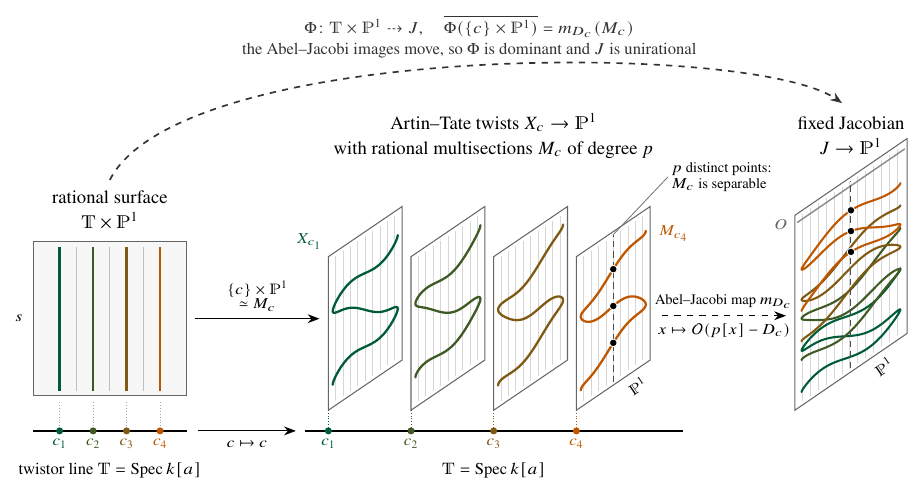}
  \caption{Sweeping out $J$ along its twistor line.
  The Abel--Jacobi images of the rational multisections $M_c$ of the Artin--Tate twists $X_c$ sweep out $J$ via a dominant rational map $\Phi\colon\TT\times\PP^1\dashrightarrow J$.
  The closed-point slices are schematic. The construction takes place over the generic point of $\TT$.
  (Figure drawn by Claude Opus~5.5.)}
  \label{fig:twistor-sweep}
\end{figure}

The main idea is to exploit the rational curves on the Artin--Tate twists of $J$, which have Artin invariant ten.
Bragg and Lieblich construct these twists from the affine line $\TT=\Spec k[a]$, the \emph{twistor line} of $J$. The group homomorphism from $\TT(k)$ onto $\Br(J)$ has finite kernel \cite[Theorem~2.2.6 and Proposition~4.3.16]{BLtw}.
Let $X_a\to\PP^1$ be the geometric generic twist, a surface over $\kab$.
Every multisection of $X_a$ has degree divisible by $p$, and the N\'eron--Severi lattice provides a smooth rational curve $M\subset X_a$ that is a multisection of degree exactly $p$ (\Cref{roots:reflection}).
Assume that $M$ is separable over the base, with generic point $D$, and that $E$ has no nonzero $p$-torsion points over $k(t)^{\sep}$.
Then the Abel--Jacobi map $x\mapsto\Oh(p[x]-D)$ from the generic fiber of $X_a$ to $E$ sends $M$ birationally onto a rational curve on $J_{\kab}$ (\Cref{sweep:abel-image}).
This curve is not defined over $k$, because only finitely many Artin--Tate torsors are splitted over a fixed separable cover of $\PP^1$, by the Lang--N\'eron theorem (\Cref{prelim:moving-covers}).
If the fibers of $J\to\PP^1$ are reduced and their component groups have order prime to $p$, then $M$ descends to a projective line over $k(a)$, the function field of the twistor line (\Cref{roots:descent}).
Its Abel--Jacobi map therefore spreads out to a dominant rational map $\TT\times\PP^1\dashrightarrow J$ from a rational surface (\Cref{fig:twistor-sweep}).
This gives a unirationality criterion for $J$ (\Cref{sweep:criterion}).

It remains to find, on every supersingular $K3$ surface of Artin invariant nine, a Jacobian elliptic fibration that satisfies the hypotheses of the criterion.
The condition on the $p$-torsion means that $E$ is geometrically supersingular or that $j(E)$ is not a $p$th power in $k(t)$ (\Cref{sweep:no-torsion}).
If $E$ is ordinary, then $M$ is automatically separable. Its degree is prime, and the geometric generic twist has no purely inseparable multisection by \Cref{prelim:no-inseparable}.
For $p\geq5$ it thus suffices to find a Jacobian elliptic fibration with $j(E)\neq0$ (\Cref{cor:Kodaira-j}).
Every Jacobian elliptic fibration has this property if $p\equiv1\pmod3$, by a result of Jang \cite{Jang} on automorphisms.
If $p\equiv2\pmod3$, we construct one with a single reducible fiber by prescribing its frame lattice, which we obtain from Kneser's neighbor method (\Cref{pencils:large}).
For $p=3$ the analogous fibration may have a supersingular generic fiber, and in that case we produce a separable rational multisection by a direct argument (\Cref{pencils:three-separable}).

\subsection*{Organization}
\Cref{sec:preliminaries} collects the preliminaries on supersingular $K3$ lattices, Artin--Tate twists and twistor lines.
\Cref{sec:separable rational multisections} deals with the separability of multisections of the geometric generic twist and their descent to the function field of the twistor line.
\Cref{sec:rational multiplication and unirationality} contains the unirationality criterion, and \Cref{sec:proof} the construction of the fibrations and the proof of \Cref{thm:main,thm:B}.

\subsection*{Conventions}
Throughout, $k$ is an algebraically closed field of characteristic $p>0$.
For a field $L$ of characteristic $p$, we write $L^{\sep}$ for a separable closure of $L$.
Cohomology with coefficients in a finite flat group scheme is fppf cohomology, and $\Br(Y)=\HH^2(Y,\GGm)$.
\subsection*{Acknowledgments}
During his graduate studies at Fudan University, the author had many discussions with Zhiyuan Li on the unirationality of supersingular $K3$ surfaces, for which he is very grateful.
We thank Christian Liedtke for his comments and for suggesting the unirationality of Enriques surfaces.

\subsection*{Declaration of AI use}
The author wrote the first draft of this paper.
Claude Code improved its organization and notation, and helped with proofreading and polishing.

The idea of the proof of \Cref{thm:main} first emerged from occasional conversations with ChatGPT~6 Astra, which produced a comparatively elementary argument with little geometric insight.
The author observed that this argument could be simplified and made more conceptual by the theory of supersingular twistor lines developed by Bragg and Lieblich \cite{BLtw}.
Building on this theory, the author then found the correct formulation of \Cref{sweep:criterion} and an alternative proof of it.
In the introduction, \Cref{fig:twistor-sweep} is generated by Claude Opus 5.5.
The same model was used to summarize the Kodaira classification of singular fibers for elliptic fibration when $\sigma =9$ (\Cref{lem:Kodaira}) from literature.
The lattice-theoretic arguments in \Cref{sec:proof} were first proposed by ChatGPT~6 Astra and ChatGPT~6.1 Sol, and were subsequently verified, corrected and substantially simplified in sustained interaction with Claude Opus~5.5.

\section{Preliminaries}\label{sec:preliminaries}
\subsection{Supersingular \texorpdfstring{$K3$}{K3} lattices and elliptic fibrations}

We first fix our terminology for (quasi-) elliptic fibrations.

\begin{definition}
   An \emph{elliptic fibration} on a smooth proper surface $X$ is a flat proper morphism $\pi\colon X \to C$ to a smooth proper curve $C$, with smooth generic fiber of genus one.
   A \emph{genus-one fibration} has the same definition with a geometrically integral generic fiber of arithmetic genus one.
   It is elliptic when the generic fiber is smooth and \emph{quasi-elliptic} when that fiber is geometrically cuspidal; the latter occurs only in characteristic $2$ and $3$.
   A genus-one fibration is \emph{Jacobian} if it admits a section $C \to X$ of $\pi$; for an elliptic fibration this means that the generic fiber is an elliptic curve.
   A multisection of degree $n$ is an integral closed subscheme $D \subset X$ such that the restriction $\pi_{D} \colon D \to C$ is finite flat of degree $n$.
\end{definition}
All fibrations considered here are relatively minimal and have connected fibers.
In characteristic $3$ we work with genus-one fibrations throughout; for $p\geq5$ every genus-one fibration is elliptic.

If $X$ is a $K3$ surface, then $C = \PP^1_k$.
The genus-one fibrations on a supersingular $K3$ surface $X$ can be described in terms of its N\'eron--Severi lattice $\NS(X)$, which is even of signature $(1,21)$.
An integral fiber of a genus-one fibration on $X$ supplies an isotropic class $f \in \NS(X)$, which we call the fiber class.
The fiber class is primitive (see \cite[Ch.~2, Remark~3.13(i)]{Huy} for example).
A class $D \in \NS(X)$ is of \emph{fiber degree} $d$ if $(D, f)=d$.
Any section of $\pi$ produces a $(-2)$-class $D$ of fiber degree one, and the pair $(D,f)$ spans a hyperbolic plane $U = \langle D +f, f \rangle$, which is a primitive sublattice of $\NS(X)$.
For a primitive isotropic fiber class $f$, the \emph{frame lattice} is $f^{\bot}/\ZZ f$.
For a Jacobian fibration it is the negative definite orthogonal complement of the copy of $U$ supplied by the fiber and section.
Roots in the frame lattice of $f$ span the root lattice of reducible fibers
\[
R_f \coloneqq \ZZ \langle r \in f^{\bot}/\ZZ f \mid r^2 = -2 \rangle,
\]
whose ADE type is called the \emph{frame root system}.

\begin{lemma}\label{lem:hyperbolic vs section}
 Let $\pi \colon X \to \PP^1$ be a genus-one fibration on a $K3$ surface $X$. Then $\pi$ is Jacobian if and only if there is a primitive hyperbolic plane $U \subset \NS(X)$ such that $U$ contains the fiber class of $\pi$.
\end{lemma}
\begin{proof}
  Let $f$ be the fiber class.
  It is primitive by \cite[Ch.~2, Remark~3.13(i)]{Huy}, applied to an integral fiber.
  If $C_0$ is a section, then $C_0^2=-2$ and $(C_0,f)=1$, so $f$ and $C_0+f$ span a hyperbolic plane $U$; it is unimodular, hence primitive.
  Conversely, let $U\subset\NS(X)$ be a hyperbolic plane containing $f$.
  Since $f$ is primitive in $\NS(X)$, it is primitive in $U$, and as $U$ is unimodular there is $D\in U$ with $(D,f)=1$.
  Let $C$ be the generic fiber of $\pi$, a regular integral proper curve over $\bk(t)$ with $\HH^0(C,\Oh_C)=\bk(t)$, arithmetic genus one and $\omega_C\cong\Oh_C$, the last by adjunction since $\omega_X\cong\Oh_X$.
  The line bundle $L=\Oh_X(D)|_C$ has degree one.
  By Serre duality, $h^1(C,L)=h^0(C,L^{-1})=0$, since $L^{-1}$ has negative degree, and Riemann--Roch gives $h^0(C,L)=1$.
  A nonzero section of $L$ cuts out an effective divisor of degree one, which is a $\bk(t)$-rational point of $C$.
  Its closure in $X$ is an integral curve mapping finitely and birationally onto the normal curve $\PP^1$, hence a section.
  For an elliptic fibration in characteristic different from $2$ and $3$, this is \cite[Ch.~11, Definition~4.3 and the discussion after it]{Huy}: $\pi$ admits a section if and only if its index is one.
\end{proof}
The proof of \Cref{lem:hyperbolic vs section} uses only that the generic fiber is a regular curve of arithmetic genus one with trivial dualizing sheaf, so it applies to quasi-elliptic fibrations in characteristic $3$.

We use the following lattice-theoretic abstraction of the supersingular case.
Let $p$ be an odd prime and $1 \leq \sigma \leq 10$ an integer.
A lattice is $p$-elementary when its discriminant group is killed by $p$.
\begin{definition}
  An even $p$-elementary lattice $N_{\sigma}$ of rank $22$, signature $(1,21)$ and discriminant $\disc N_{\sigma} = - p^{2\sigma}$ is called a \emph{supersingular $K3$ lattice} of Artin invariant $\sigma$ in characteristic $p$.
\end{definition}
As $N_{\sigma}$ is $p$-elementary, the quadratic form induced on the $\FF_p$-vector space $N_{\sigma}/p N_{\sigma}$ has rank $22 - 2 \sigma$, and its radical is $pN_\sigma^\vee/pN_\sigma$.
Rudakov and Shafarevich prove the $p$-elementarity of the N\'eron--Severi lattice and the uniqueness of the lattice $N_\sigma$ in \cite[\S8 and \S1, respectively]{RS}. See also \cite[Ch.~17, Propositions~2.19--2.20 and Corollary~2.21]{Huy}.
\begin{theorem}[Rudakov--Shafarevich]\label{thm:RS}
Let $p$ be an odd prime and let $1\leq\sigma\leq10$.
\begin{enumerate}[label=\textup{(\roman*)},leftmargin=2.2em]
  \item If $X$ is a supersingular $K3$ surface with $\sigma(X)=\sigma$ over an algebraically closed field of characteristic $p$, then $\NS(X)$ is a supersingular $K3$ lattice of Artin invariant $\sigma$ in characteristic $p$.
  \item Any two supersingular $K3$ lattices of Artin invariant $\sigma$ in characteristic $p$ are isometric.
\end{enumerate}
In particular, for $X$ as in \textup{(i)}, every even lattice of signature $(1,21)$ with discriminant group $(\ZZ/p)^{2\sigma}$ is isometric to $\NS(X)$.
\end{theorem}
The last assertion is the form in which the theorem is used below: such a lattice has discriminant $-p^{2\sigma}$, so it is a supersingular $K3$ lattice of Artin invariant $\sigma$.
Moreover, for any $N_{\sigma}$, there is a supersingular $K3$ surface $X$ over an algebraically closed field $k$ with $\cha(k)=p$ and $\NS(X) \cong N_{\sigma}$.

Any elliptic fibration of a supersingular $K3$ surface of Artin invariant $10$ admits no sections by Ekedahl and van der Geer \cite[Proposition 13.1]{EvdG11}; see \cite[Corollary 1.6]{KS} for another proof.
The same holds for genus-one fibrations by \Cref{lem:hyperbolic vs section}, since $N_{10}$ contains no hyperbolic plane (see the proof of \Cref{prelim:jacobian}).
For $\sigma \leq 9$, in contrast, Jacobian genus-one fibrations exist.
For $p\geq 5$ this is \cite[Proposition~3.9]{Lie}, and the genus-one formulation in characteristic $3$ is stated in \cite[Remark~3.11]{Lie}.
\begin{proposition}\label{prelim:jacobian}
In odd characteristic a supersingular $K3$ surface of Artin invariant at most $9$ has a Jacobian genus-one fibration. On a supersingular $K3$ surface of Artin invariant $10$, no genus-one fibration admits a section.
\end{proposition}
\begin{proof}  
Let $X$ be a supersingular $K3$ surface with $\sigma=\sigma(X)$, and identify $\NS(X)$ with $N_\sigma$ by \Cref{thm:RS}.
For $\sigma\leq9$ the quadratic form on $N_{\sigma}/pN_{\sigma}$ has rank at least four.
A nondegenerate quadratic space of dimension at least three over $\Fp$ has a nonzero isotropic vector, so there is an isotropic vector in $N_\sigma/pN_\sigma$ outside the radical.
The derivative of the quadratic form is nonzero there, so Hensel's lemma lifts it to an integral isotropic $p$-adic vector whose pairings with $N_\sigma\otimes\ZZ_p$ generate $\ZZ_p$.
The rational quadratic space is isotropic by Meyer's theorem.
Its projective isotropic quadric is rational and satisfies weak approximation.
A rational point close to the chosen $p$-adic line therefore gives a primitive integral isotropic vector $f$ whose divisibility is prime to $p$.
The divisibility of a primitive vector divides the exponent $p$ of the discriminant group, so it is one.
Choosing $g$ with $(f,g)=1$ and replacing $g$ by $g-(g^2/2)f$ generates a hyperbolic plane $U$.
Choose an ample class $H$.
Replacing $f,g$ simultaneously by $-f,-g$ if necessary makes $(f,H)>0$, and Riemann--Roch then makes $f$ effective.
Applying the same root reflections to both vectors puts $f$ in the nef cone and preserves $(f,g)=1$.
The primitive nef isotropic class $f$ defines a genus-one fibration.
Its fiber class lies in the hyperbolic plane spanned by $f$ and $g$, so the fibration is Jacobian by \Cref{lem:hyperbolic vs section}.
It may be quasi-elliptic, which is allowed in characteristic three.

Suppose now that $N_{10}$ contains a hyperbolic plane $U$.
It is unimodular, hence an orthogonal summand, and its complement $L$, a negative definite even lattice of rank twenty, would have discriminant group $(\ZZ/p)^{20}$.
Its Gram matrix would be divisible by $p$, and $L(1/p)$ would be an even unimodular definite lattice of rank twenty.
Such a lattice does not exist, since its signature would not be divisible by eight.
By \Cref{lem:hyperbolic vs section}, no genus-one fibration on a surface of Artin invariant $10$ has a section.
\end{proof}

An analysis of the lattice $N_{10}$ yields rational multisections on genus-one fibrations in Artin invariant $10$.
This seems to be well known to experts; we record it for the convenience of the reader.

\begin{lemma}\label{roots:lattice}
  Put $N=N_{10}$.
  For any primitive isotropic class $f \in N$, there is another isotropic class $g \in N$ such that
\[
  N=U(p)\perp L,\qquad U(p)=\ZZ f+\ZZ g,
\]
where the modulo $p$ quadratic form on $L/pL$ has rank two.
Moreover, for $v\in L$ with $v^2\equiv-2\pmod p$, the class
\begin{equation}\label{eq:root}
    r_v=g+v+\frac{-2-v^2}{2p}f
\end{equation}
is a root of fiber degree $p$, that is, $r_v^2=-2$ and $(r_v,f)=p$.
Such $v$ exist, and every root $r\in N$ with $(r,f)=p$ is $r_v$ for $v$ the component of $r$ in $L$.
\end{lemma}
\begin{proof}
Since $f$ is primitive, its divisibility divides the exponent $p$ of the discriminant group.
If the divisibility of $f$ is equal to one, then, as $f$ is isotropic, $N$ contains a hyperbolic plane, which contradicts the proof of \Cref{prelim:jacobian}.
Thus $(f,N)=p\ZZ$.

Choose $w\in N^\vee$ with $(w,f)=1$.
Since $pN^\vee\subset N$, the square of $pw$ is even and divisible by $p$.
Thus
\[
            g=pw-\frac{(pw)^2}{2p}f
\]
is integral, isotropic, has $(f,g)=p$, and pairs with $N$ in $p\ZZ$.
Orthogonal projection along $f,g$ is therefore integral, which gives the decomposition $N=U(p)\perp L$.
The form on $U(p)$ vanishes modulo $p$, and the form on $N/pN$ has rank $22-2\cdot10=2$, so the form on $L/pL$ has rank two.

Every nondegenerate binary quadratic form over $\Fp$ represents each nonzero scalar: this is immediate for the split form and follows from the norm map $\FF_{p^2}^\times\to\Fp^\times$ for the nonsplit form.
Hence there is $v\in L$ with $v^2\equiv-2\pmod p$.
Since $v^2$ is even and $p$ is odd, $(-2-v^2)/2p$ is an integer, and a direct computation shows that the class $r_v$ of \eqref{eq:root} satisfies $r_v^2=-2$ and $(r_v,f)=p$.
Conversely, write a root $r$ with $(r,f)=p$ as $r=xf+yg+v$ with $v\in L$.
Then $py=(r,f)=p$, so $y=1$, and $-2=r^2=2px+v^2$ determines $x$.
\end{proof}

\begin{proposition}\label{roots:reflection}
Let $X$ be a supersingular $K3$ surface with $\sigma(X)=10$ in odd characteristic $p$, and let $\pi\colon X\to\PP^1$ be a genus-one fibration with fiber class $f\in N=\NS(X)$.
For a $(-2)$-curve $C$, write $s_C(x)=x+(x,C)C$ for the reflection in $C$.
\begin{enumerate}[label=\textup{(\roman*)},leftmargin=2.2em]
  \item Every multisection of $\pi$ has degree divisible by $p$.
  \item There exist roots $r\in N$ with $(r,f)=p$, and every such root is effective.
  \item For every root $r\in N$ with $(r,f)=p$ there are vertical $(-2)$-curves $C_1,\dots,C_k$ such that $s_{C_k}\circ\dots\circ s_{C_1}(r)$ is the class of a smooth rational curve $M\subset X$.
    The curve $M$ is a multisection of $\pi$ of degree $p$.
\end{enumerate}
\end{proposition}
\begin{proof}
The fiber class $f$ is primitive and isotropic, so \Cref{roots:lattice} applies.

(i) The degree of a multisection $M$ is $(M,f)\in(f,N)=p\ZZ$.

(ii) Roots with $(r,f)=p$ exist by \Cref{roots:lattice}.
For such a root, Riemann--Roch gives $h^0(r)+h^0(-r)\geq\chi(\Oh_X(r))=2+r^2/2=1$ \cite[Chapter~2, Section~1.4]{Huy}, so $r$ or $-r$ is effective.
Since $f$ is nef and $(-r,f)=-p<0$, the class $-r$ is not effective.

(iii) Fix an ample class $H$.
If $C$ is a vertical $(-2)$-curve with $(r,C)<0$, then $s_C(r)$ is again a root with $(s_C(r),f)=p$, because $(C,f)=0$.
It is effective by \textup{(ii)}, and $(s_C(r),H)=(r,H)+(r,C)(C,H)<(r,H)$.
Degrees of nonzero effective classes with respect to $H$ are positive integers, so after finitely many such reflections, as in \cite[Chapter~8, Remark~2.13]{Huy}, we reach a root $r'$ with $(r',f)=p$ and $(r',C)\geq0$ for every vertical $(-2)$-curve $C$.
Every vertical integral curve is either a $(-2)$-curve or an integral fiber, and the latter meets $r'$ in $p$.
Hence $(r',C)\geq0$ for every vertical curve $C$.

Write the effective class $r'$ as $M+D$, where $D$ is effective and vertical and $M$ is effective with no vertical components.
Then $(M,f)=p$.
Every horizontal integral curve has positive degree in $(f,N)=p\ZZ$ over $\PP^1$, so $M$ is a single integral curve with multiplicity one, and $M^2\geq-2$ \cite[Chapter~2, Section~1.3]{Huy}.
On the other hand,
\[
  -2-M^2=r'^2-M^2=(r'+M,D)=(r',D)+(M,D)\geq0,
\]
since $r'$ is nonnegative on vertical curves and $M$ has no vertical components.
Hence $M^2=-2$, so $M\simeq\PP^1$, and $(r',D)=(M,D)=0$, whence $D^2=(r',D)-(M,D)=0$.
Write $D=\sum_vD_v$ with $D_v$ supported on the fiber over $v$.
By Zariski's lemma \cite[Chapter~11, Proposition~1.6(iv)]{Huy}, whose proof does not use the restriction on the characteristic, $D_v^2\leq0$ with equality only if $D_v$ is a rational multiple of the fiber.
As $D^2=\sum_vD_v^2=0$, every $D_v$ is a nonnegative rational multiple $q_vf$, and $0=(M,D)=p\sum_vq_v$ forces $D=0$.
Thus $r'=[M]$.
Finally, $M$ is integral and not vertical, so $M\to\PP^1$ is finite and flat over the regular curve $\PP^1$, of degree $(M,f)=p$.
\end{proof}
\subsection{Singular fibers and \texorpdfstring{$j$}{j}-invariants in Artin invariant \texorpdfstring{$9$}{9}}\label{subsec:Kodaira}
Throughout this subsection, $p$ is odd, $J$ is a supersingular $K3$ surface with $\sigma(J)=9$, and $\pi\colon J\to\PP^1$ is a Jacobian elliptic fibration with fiber class $f$ and generic elliptic curve $E/\bk(t)$.
For a fiber $F_v$, let $m_v$ be its number of components, let $e(F_v)$ be its Euler number, and let $R_v$ be the root lattice spanned by its components that do not meet the zero section.
Write $\Phi_v$ for the component group of the N\'eron model of $E$ at $v$, $\Delta$ for the discriminant of a minimal Weierstrass model at $v$, and $j=j(E)\in\bk(t)$ for the $j$-invariant.
\begin{lemma}\label{lem:Kodaira}
  The frame root system is $R_f=\bigoplus_vR_v$.
  It is an orthogonal sum of at most two lattices, each isomorphic to $A_1$ or $A_2$, and if $p\geq5$, it is one of $0$, $A_1$, $A_1^{\oplus2}$, $A_2$.
  Every singular fiber $F_v$ of $\pi$ is of one of the following six types, with the listed invariants.
  \begin{center}\normalfont
    \renewcommand{\arraystretch}{1.2}
    \begin{tabular}{cccccl}
      \hline
      $F_v$ & $R_v$ & $m_v$ & $\Phi_v$ & $e(F_v)$ & behavior of $j$ at $v$ \\
      \hline
      $\KI_1$ & $0$ & $1$ & $0$ & $1$ & pole of order $1$ \\
      $\KI_2$ & $A_1$ & $2$ & $\ZZ/2$ & $2$ & pole of order $2$ \\
      $\KI_3$ & $A_2$ & $3$ & $\ZZ/3$ & $3$ & pole of order $3$ \\
      $\KII$ & $0$ & $1$ & $0$ & $2$ & regular, $j(v)=0$ \\
      $\KIII$ & $A_1$ & $2$ & $\ZZ/2$ & $3$ & regular, $j(v)=1728$ \\
      $\KIV$ & $A_2$ & $3$ & $\ZZ/3$ & $4$ & regular, $j(v)=0$ \\
      \hline
    \end{tabular}
  \end{center}
  The values $j(v)$ for the types $\KII$, $\KIII$ and $\KIV$ are claimed only for $p\geq 5$.
  In particular, all fibers are reduced, $|\Phi_v|\leq3$, and $|\Phi_v| = 3$ only if $R_v\cong A_2$.
  If $p\geq5$, then $e(F_v)=\ord_v\Delta$ and $\sum_ve(F_v)=24$.
\end{lemma}
\begin{proof}
  An elliptic $K3$ surface has a singular fiber by \cite[Ch.~11, Remark~1.5(ii)]{Huy}.
  Since $\pi$ has a section, it has no multiple fibers, so every fiber is of Kodaira type \cite[Section~4.1]{SS}.
  By Riemann--Roch, a root $r\in f^\perp$ or its negative is effective, and its components are vertical because $f$ is nef.
  Hence $R_f=\bigoplus_vR_v$; see \cite[Lemma~8.3]{SS}.

  Let $W=f^\perp/\ZZ f$ be the frame lattice.
  The fiber class and the zero section span a hyperbolic plane in $\NS(J)$ with orthogonal complement isometric to $W$, so $W$ has rank $20$ and discriminant group $(\ZZ/p)^{18}$ by \Cref{thm:RS}.
  Thus the form on $W/pW$ has rank two, and every Gram matrix of vectors in $W$ has rank at most two modulo $p$.
  A connected Dynkin diagram with at least three vertices contains $A_3$, whose determinant $-4$ is a unit modulo $p$.
  Thus the summands of $R_f$ are $A_1$ and $A_2$.
  Their Gram matrices are nonzero modulo $p$, so there are at most two summands.
  If $p\geq5$, their determinants are prime to $p$, so $\rank R_f\leq2$.
  Kodaira's list leaves the six fiber types above, and their component groups are listed in \cite[Lemma~7.3]{SS}.

  Since $p\neq2$, the function $j$ has a pole at $v$ exactly when $E$ has potentially multiplicative reduction at $v$, that is, when $F_v$ is of type $\KI_n$ or $\KI_n^*$ with $n\geq1$, and the pole then has order $n$ \cite[Section~4]{SS}.
  This gives the poles in the first three rows and the regularity in the last three.

  Let now $p\geq5$.
  The values of $j$ follow from a minimal Weierstrass form $y^2=x^3+a_4x+a_6$, where $j=1728\cdot4a_4^3/(4a_4^3+27a_6^2)$, and the orders of $a_4,a_6,\Delta$ listed in \cite[Section~4.5, Table~1]{SS}.
  For instance, at a fiber of type $\KIII$ one has $\ord_va_4=1$, $\ord_va_6\geq2$ and $\ord_v\Delta=3$, so $\ord_v(j-1728)=2\ord_va_6-3>0$.
  There is no wild ramification for $p\geq5$ \cite[Section~4.6]{SS}.
  Hence $e(F_v)=\ord_v\Delta$, and \cite[Theorem~6.10]{SS} gives $\sum_ve(F_v)=e(J)=24$.
\end{proof}
\begin{remark}
The first step of this proof uses only Riemann--Roch and the nefness of $f$, so $R_f=\bigoplus_vR_v$ holds for every genus-one fibration on a $K3$ surface.
\end{remark}

\begin{corollary}\label{cor:Kodaira-j}
  \leavevmode
  \begin{enumerate}[label=\textup{(\roman*)},leftmargin=2.2em]
    \item If $j$ is nonconstant, then its poles are the fibers of type $\KI_n$, of order $n\leq3$.
      Hence $j\notin\bk(t)^p$ if $p\geq5$, and also if $p=3$ and no fiber is of type $\KI_3$.
    \item If $j$ is constant, then $j=0$, and no fiber is of type $\KI_n$ with $n\geq1$.
      If moreover $p\geq5$, the singular fibers are either twelve fibers of type $\KII$, or ten fibers of type $\KII$ and one of type $\KIV$.
  \end{enumerate}
  In particular, if $R_f$ is $A_1$ or $A_1^{\oplus2}$, then either $j\notin\bk(t)^p$, or $p=3$ and $j=0$.
\end{corollary}
\begin{proof}
  (i) This follows from \Cref{lem:Kodaira}, since a nonconstant function has a pole, and a $p$th power in $\bk(t)$ has pole orders divisible by $p$.

  (ii) A constant has no poles, so no fiber is of type $\KI_n$ with $n\geq1$.
  As in the proof of \Cref{lem:Kodaira}, at least one fiber is singular.

  Let $p\geq5$.
  Fibers of types $\KII$ and $\KIV$ force $j=0$, and fibers of type $\KIII$ force $j=1728\neq0$.
  If $j=1728$, there are $24/3=8$ fibers of type $\KIII$, contradicting $\rank R_f\leq2$.
  Hence $j=0$, and there are at most one fiber of type $\KIV$ and no fibers of type $\KIII$.
  The equation $2\#\KII+4\#\KIV=24$ gives the two cases.

  Let $p=3$, and suppose that $j\neq0$.
  As $j\neq0=1728$, the automorphism group of the geometric generic fiber is $\{\pm1\}$, so $E$ is a quadratic twist of a constant elliptic curve, and every singular fiber of $\pi$ is of type $\KI_0^*$ \cite[Section~5.4]{SS}.
  This type does not occur by \Cref{lem:Kodaira}, a contradiction.

  For the last assertion, let $R_f$ be $A_1$ or $A_1^{\oplus2}$, so that no fiber is of type $\KI_3$ or $\KIV$.
  If $j$ is nonconstant, then $j\notin\bk(t)^p$ by \textup{(i)}.
  If $j$ is constant, then $j=0$ and $p=3$ by \textup{(ii)}, since for $p\geq5$ a constant $j$ allows no fiber of type $\KI_2$ or $\KIII$.
\end{proof}

\subsection{Artin--Tate twistor families}\label{sec:artin-tate}
We recall the geometric inputs from the theory of Artin--Tate twistors established in \cite{BLtw}.

Throughout this subsection $p=\cha(\bk)$ is odd, and $\pi\colon J\to\PP^1_\bk$ is a supersingular Jacobian elliptic $K3$ surface with generic fiber $E/\bk(t)$ and Artin invariant $\sigma=\sigma(J)$.
The letter $\Omega$ denotes an algebraically closed field containing $\bk$, possibly $\Omega=\bk$.
We write $J_\Omega$ and $\pi_\Omega$ for the base changes and $E_{\Omega(t)}=E\otimes_{\bk(t)}\Omega(t)$.
The representability of Picard functor provides the geometric invariance of Artin invariants.
\begin{lemma}\label{prelim:base-change}
Let $X$ be a $K3$ surface over $\bk$, with structure morphism $q\colon X\to\Spec\bk$.
The Picard scheme $\Pic_{X/\bk}$ is the constant group scheme on the torsion-free group $\Pic(X)\cong \NS(X)$, and $\Pic(X)\to\Pic(X_\Omega)$ is an isometry.
In particular, $X_\Omega$ is supersingular if and only if $X$ is, and then $\sigma(X_\Omega)=\sigma(X)$.
Moreover, $q_*\mup=\mup$, and the fppf sheaf $R^1q_*\mup$ vanishes.
\end{lemma}
\begin{proof}
The Picard scheme is locally of finite type with tangent space $\HH^1(X,\Oh_X)=0$, so it is \'etale over $\bk$.
Since $\bk$ is algebraically closed, it is a disjoint union of copies of $\Spec\bk$.
As $X$ has a $\bk$-point, its $\bk$-points form $\Pic(X)$ and its $\Omega$-points form $\Pic(X_\Omega)$, so it is the constant group scheme on $\Pic(X)$ and $\Pic(X)\to\Pic(X_\Omega)$ is bijective.
The Picard group of a $K3$ surface is torsion-free and equals its N\'eron--Severi group, and intersection numbers do not change under field extension.
Since $q_*\GGm=\GGm$, the Kummer sequence (in fppf topology) \[
0\longrightarrow \mup\longrightarrow \GGm\longrightarrow \GGm\longrightarrow 0 \] gives $q_*\mup=\mup$ and identifies $R^1q_*\mup$ with the $p$-torsion of $\Pic_{X/\bk}$, which is zero.
\end{proof}

Let $l$ be an algebraically closed field and $A$ an elliptic curve over $l(t)$.
The \emph{Tate--Shafarevich group} of $A$ is
\[
  \Sha(A) = \ker\left( \HH^1(l(t), A) \longrightarrow \prod_{v \in \PP^1(l)} \HH^1(l(t)_v, A) \right),
\]
where $l(t)_{v}$ is the completion of $l(t)$ at the $l$-point $v$.
The \emph{period} of a torsor $T$ is the order of its class $[T]$ in $\HH^1(l(t),A)$, and its \emph{index} is the positive generator of the degrees of its closed points.

\begin{lemma}\label{lem:ArtinTate}
Let $X$ be a supersingular $K3$ surface over an algebraically closed field $K$ of characteristic $p$.
Then $\Br(X)$ is killed by $p$, and the Kummer sequence induces a short exact sequence of groups
\begin{equation}\label{eq:unipotent Brauer}
  0 \longrightarrow \frac{p\NS(X)^{\vee}}{p\NS(X)} \longrightarrow \GGa(K) \longrightarrow \Br(X) \longrightarrow 0.
\end{equation}
In particular, $\Br(X)\neq0$.
If $X=J_\Omega$, there is moreover a canonical isomorphism $\Br(J_\Omega) \simeq \Sha(E_{\Omega(t)})$.
\end{lemma}
\begin{proof}
The group $\Br(X)$ is killed by $p$ by \cite[Theorem~2.2.6]{BLtw}. The top row of its Kummer diagram gives \eqref{eq:unipotent Brauer}, with $U^2(X,\mup)$ identified with $\GGa(K)$.
In the elliptic case these results go back to \cite[Theorem~(4.3)(ii) and Corollary~(4.4)]{Artin}, under his flat duality hypothesis. Milne established the duality consequences needed here \cite[Remark~5.3]{Milne}.
The group on the left of \eqref{eq:unipotent Brauer} is finite and $K$ is infinite, so $\Br(X)\neq0$.
The isomorphism $\Br(J_\Omega) \simeq \Sha(E_{\Omega(t)})$ is the theorem of Artin--Tate and Grothendieck \cite[Theorem~3.1]{AT}; see also \cite[Theorem~4.3.8]{BLtw}.
By \Cref{prelim:base-change}, $J_\Omega$ is a supersingular Jacobian elliptic $K3$ surface over $\Omega$ with generic fiber $E_{\Omega(t)}$, so this applies over $\Omega$.
\end{proof}

\begin{proposition}\label{prelim:artin-tate}
For $\alpha\in\Br(J_\Omega)$, let $X(\alpha)\to\PP^1_\Omega$ be the relatively minimal regular model of the $E_{\Omega(t)}$-torsor whose class corresponds to $\alpha$ under \Cref{lem:ArtinTate}.
Then $X(\alpha)\to\PP^1_\Omega$ is \'etale locally on $\PP^1_\Omega$ isomorphic to $\pi_\Omega$, and $X(\alpha)$ is a supersingular $K3$ surface over $\Omega$.
If $\alpha\neq0$, it has Artin invariant $\sigma+1$.
If $\alpha=0$, then $X(\alpha)\cong J_\Omega$ over $\PP^1_\Omega$.
\end{proposition}
\begin{proof}
By \Cref{prelim:base-change}, $J_\Omega$ is a supersingular Jacobian elliptic $K3$ surface over $\Omega$ with Artin invariant $\sigma$ and generic fiber $E_{\Omega(t)}$.
Replacing $\bk$ by $\Omega$, we may assume $\Omega=\bk$.

If $\alpha=0$, the torsor is trivial, hence isomorphic to $E$.
The surface $J$ is a regular model of $E$, and it is minimal since a $K3$ surface contains no $(-1)$-curves.
By uniqueness of relatively minimal regular models, $X(0)\cong J$ over $\PP^1_\bk$.

For $\alpha\neq0$, choose a lift of $\alpha$ to $U^2(J,\mup)$ by \cite[Theorem~2.2.6]{BLtw}. The corresponding $\mup$-gerbe $\cG_\alpha$ deforms the trivial gerbe.
By \cite[Definition~4.3.2(2), Lemma~4.3.5 and Proposition~4.3.16]{BLtw}, the surface $X(\alpha)$ is the coarse moduli space of stable $\cG_\alpha$-twisted sheaves of Mukai vector $v=(0,f,1)$, where $f$ is the fiber class and Gieseker-stability is taken with respect to a $v$-generic polarization.
In particular, $X(\alpha)$ is \'etale locally isomorphic to $J$ over $\PP^1_\bk$, and it is a supersingular $K3$ surface by \cite[Theorem~4.1.19(1)]{BLtw}.

Let $S\simeq\PP^1_\bk$ be the zero section of $J$. Since $\Br(S)=0$, there is a twisted line bundle on $\cG_\alpha|_S$.
Its pushforward to $\cG_\alpha$ has Mukai vector $w=(0,[S],s)$ for some $s$, so $v\cdot w=(f,S)=1$. By \cite[Theorem~4.1.19(3)]{BLtw}, the moduli gerbe is trivial.
The universal sheaf therefore gives a Fourier--Mukai equivalence $D(J,\alpha)\simeq D(X(\alpha))$ by \cite[Theorem~4.1.19(2)]{BLtw}.
Derived equivalent twisted supersingular $K3$ surfaces have the same twisted Artin invariant \cite[Corollary~3.4.3]{Bragg}.
A nonzero Brauer class raises the Artin invariant by one \cite[Corollary~3.4.23]{BLtw}, \cite[Proposition~3.2.7]{Bragg}.
Hence $\sigma(X(\alpha))=\sigma_0(J,\alpha)=\sigma+1$.
\end{proof}

Here ``\'etale locally'' cannot be replaced by a finite \'etale cover, since $\PP^1_\Omega$ is simply connected.
We call the surfaces $X(\alpha)$ with $\alpha\neq0$ the \emph{Artin--Tate twists} of $J$ over $\Omega$, or simply of $J$ when $\Omega=\bk$.
The generic fiber of $X(\alpha)$ is the $E_{\Omega(t)}$-torsor with class $\alpha$.
We reserve the word torsor for such genus-one curves over function fields, and for torsors under group schemes, and never use it for the surfaces.

Let $T$ be the generic fiber of $X(\alpha)$.
By \Cref{lem:ArtinTate}, $p[T]=0$, so the pushout $[p]_*[T]$ of $T$ along the multiplication $[p]\colon E_{\Omega(t)}\to E_{\Omega(t)}$ is trivial.
An identification of $[p]_*T=T\times^{E_{\Omega(t)},[p]}E_{\Omega(t)}$ with the trivial torsor gives a morphism $T\to E_{\Omega(t)}$ over $\Omega(t)$ inducing $[p]$ on Jacobians.
It is geometrically a translate of $[p]$, so it is finite of degree $p^2$, and it is unique up to translation by $E_{\Omega(t)}(\Omega(t))$.
We call it a \emph{$[p]$-map} on $T$ and still denote it by $[p]$.

\begin{proposition}\label{prop:AT twistor family}
\leavevmode
\begin{enumerate}[label=\textup{(\roman*)},leftmargin=2.2em]
\item For every Artin--Tate twist $X=X(\alpha)$ of $J$ over $\Omega$, a $[p]$-map $T\to E_{\Omega(t)}$ on its generic fiber extends to a dominant rational map $X\dashrightarrow J_\Omega$ over $\PP^1_\Omega$ of degree $p^2$.
  In particular, if $X$ is unirational, then so is $J_\Omega$.
\item Let $1\le s\le9$.
  If every supersingular $K3$ surface over $\bk$ of Artin invariant $s+1$ is unirational, then so is every supersingular $K3$ surface over $\bk$ of Artin invariant $s$.
\end{enumerate}
\end{proposition}
\begin{proof}
(i) The $[p]$-map is a finite morphism $T\to E_{\Omega(t)}$ of degree $p^2$ between the generic fibers of $X$ and $J_\Omega$.
Since $X$ and $J_\Omega$ are integral and proper over $\PP^1_\Omega$ with these generic fibers, it extends to a rational map $X\dashrightarrow J_\Omega$ over $\PP^1_\Omega$, which is dominant of degree $p^2$ because it is so on generic fibers.
Composing a dominant rational map $\PP^2_\Omega\dashrightarrow X$ with it proves the second assertion.

(ii) Let $Y$ be supersingular of Artin invariant $s\le9$.
By \Cref{prelim:jacobian} it has a Jacobian genus-one fibration.
If this fibration is quasi-elliptic, which can only happen for $p=3$, then $Y$ is unirational directly as observed in \cite[\S2, Proposition~2]{RS}. After a purely inseparable base change of $\PP^1$, the normalization of the generic cuspidal curve becomes a smooth rational curve over the rational field $\bk(t^{1/p^n})$, and the section supplies a rational point on it.
Thus it is a projective line over $\bk(t^{1/p^n})$, and $Y$ is dominated by a rational surface.
If the fibration is elliptic, we may apply the preceding results with $J=Y$ and $\Omega=\bk$.
By \Cref{lem:ArtinTate}, $\Br(Y)\neq0$.
For $\alpha\neq0$ in $\Br(Y)$, the Artin--Tate twist $X(\alpha)$ has Artin invariant $s+1$ by \Cref{prelim:artin-tate}.
It is unirational by hypothesis, and $Y$ is unirational by \textup{(i)}.
\end{proof}

\begin{corollary}\label{prelim:descent-ten}
If all supersingular $K3$ surfaces over $\bk$ of Artin invariant $10$ are unirational, then all supersingular $K3$ surfaces over $\bk$ are unirational.
\end{corollary}

When the generic fiber is geometrically supersingular, the $[p]$-map is purely inseparable, and there is also a map in the other direction.
\begin{lemma}\label{lem:purely inseparable sandwich}
  Assume that the geometric generic fiber of $\pi \colon J \to \PP^1_k$ is supersingular. Then any Artin--Tate twist $X$ of $J$ over $\Omega$ admits a sequence of purely inseparable dominant rational $\Omega$-maps, each of degree $p^2$:
  \[
  J^{(1/p^2)}_{\Omega} \dashrightarrow X \dashrightarrow J_{\Omega}.
  \]
  Here, the inverse Frobenius twist $J^{(1/p^2)}$ is the base change of $J$ along the inverse of the $p^2$-power Frobenius automorphism of the perfect field $k$.
\end{lemma}
\begin{proof}
By \Cref{prop:AT twistor family}\textup{(i)}, a $[p]$-map on the generic fiber of $X$ induces a dominant rational $\Omega$-map $f\colon X\dashrightarrow J_\Omega$ of degree $p^2$.
After extending $\Omega(t)$ to an algebraic closure and trivializing the generic torsor, this map is a translate of multiplication by $p$ on a supersingular elliptic curve.
Its geometric fibers therefore have a single point, since a supersingular elliptic curve has no nonzero geometric $p$-torsion points.
Thus the generic-fiber map, and hence $f$, is purely inseparable.

Identify $K=\Omega(J_\Omega)$ with its image in $L=\Omega(X)$ under $f^*$.
The extension $L/K$ is purely inseparable of degree $p^2$.
For every $u\in L$, its minimal polynomial over $K$ has the form $z^{p^r}-u^{p^r}$, with $p^r=[K(u):K]\leq p^2$, so $u^{p^2}\in K$.
Consequently, inside a fixed algebraic closure of $K$,
\[
  K\subset L\subset K^{1/p^2}.
\]
Since $\Omega$ is perfect, the function field of $(J_\Omega)^{(1/p^2)}$ is naturally $K^{1/p^2}$ as an $\Omega$-field.
Moreover, inverse Frobenius twisting commutes with the extension $k\subset\Omega$, so $(J_\Omega)^{(1/p^2)}\cong J^{(1/p^2)}_\Omega$.
The inclusion $L\subset K^{1/p^2}$ therefore defines a purely inseparable dominant rational $\Omega$-map $g\colon J^{(1/p^2)}_\Omega\dashrightarrow X$.
As $K$ is the function field of a smooth surface over the perfect field $\Omega$, two iterations of Frobenius have degree $p^4$ \cite[Tag~0CCF]{stacks-project}.
Hence
\[
  \deg g=[K^{1/p^2}:L]
  =\frac{[K^{1/p^2}:K]}{[L:K]}
  =\frac{p^4}{p^2}=p^2,
\]
as required.
\end{proof}
\subsection{Supersingular twistor lines}
The Brauer group of a supersingular $K3$ surface $X$ is parametrized by an affine line \cite[Theorem~2.2.6]{BLtw}.
\begin{proposition}[The twistor line]\label{prelim:twistor}
Let $X$ be a supersingular $K3$ surface over $\bk$, with structure morphism $q\colon X\to\Spec\bk$.
\begin{enumerate}[label=\textup{(\roman*)},leftmargin=2.2em]
\item The fppf sheaf $R^2q_*\mup$ on the category of $\bk$-schemes is represented by a smooth commutative group scheme $\bH_X$ over $\bk$, whose identity component is isomorphic to $\GGa$.
  Fix a group isomorphism
\[
\TT=\bH_X^0\simeq\Spec\bk[a]=\A^1_\bk.
\]
We call $\TT$ the \emph{twistor line} of $X$ and $a$ its \emph{coordinate}; two coordinates differ by a scalar in $\bk^\times$.
\item There is a $\mup$-gerbe
\[
              \cG\longrightarrow X\times\TT,
\]
unique up to equivalence, whose class in $\HH^2(X\times\TT,\mup)$ maps to the tautological section $\TT\hookrightarrow\bH_X$ of $R^2q_*\mup$ over $\TT$.
For a field $K\supseteq\bk$ and a point $c\in\TT(K)=K$, let $\cG_c$ be the pullback of $\cG$ along $(\mathrm{id}_X,c)\colon X_K\to X\times\TT$, and let $\alpha_c\in\Br(X_K)$ be the image of its class under the Kummer map $\HH^2(X_K,\mup)\to\Br(X_K)$.
For a $\bk$-embedding $\iota\colon K\to K'$, the class $\alpha_c$ pulls back to $\alpha_{\iota(c)}$ along $X_{K'}\to X_K$.
\item For every algebraically closed $\Omega\supseteq\bk$, the map $\Omega\to\Br(X_\Omega)$, $c\mapsto\alpha_c$, is a surjective homomorphism, equivariant for $\Aut(\Omega/\bk)$.
  Its kernel is a finite subgroup $Z_X\subset\bk$, the same for every $\Omega$.
  In particular, $\alpha_c\neq0$ for every $c\in\Omega$ transcendental over $\bk$.
\end{enumerate}
\end{proposition}
\begin{proof}
\textup{(i)} The smooth group scheme representing $R^2q_*\mup$ and the identification $\bH_X^0\simeq\GGa$ are given by \cite[Proposition~2.2.4]{BLtw}. See also \cite[Theorem~4.1]{LieTw}.
A group automorphism of $\GGa$ is an automorphism of the affine line fixing the origin, so it is multiplication by a scalar.

\textup{(ii)} 
Lieblich lifts the inclusion of the identity component along the Leray edge map for an elliptic $K3$ surface of Artin invariant $10$ \cite[Section~6, p.~19]{LieTw}. The construction for a general supersingular $K3$ surface is given after Lemma~2.2.5 in \cite{BLtw}. We recall the argument here.
We use the Leray spectral sequence $E_2^{i,j}=\HH^i(\TT,R^j\pr_{2*}\mup)\Rightarrow \HH^{i+j}(X\times\TT,\mup)$ of the projection $\pr_2\colon X\times\TT\to\TT$.
Here $R^j\pr_{2*}\mup$ is the restriction of $R^jq_*\mup$ to $\TT$.
By \Cref{prelim:base-change}, $R^1q_*\mup=0$, and $q_*\mup=\mup$.
The line $\TT$ is affine, so $\HH^3(\TT,\mup)=0$ by \cite[Lemma~2.2.5]{BLtw}.
The Kummer sequence on $\TT$ gives $\HH^2(\TT,\mup)=0$, since $\Pic(\TT)=0$ and $\Br(\TT)\subset\Br(\bk(a))=0$ by Tsen's theorem.
Thus $E_2^{i,j}=0$ whenever $i+j\in\{2,3\}$ and $j\leq1$, and the edge map
\[
 \HH^2(X\times\TT,\mup)\longrightarrow
   \HH^0(\TT,R^2\pr_{2*}\mup)=\Hom_\bk(\TT,\bH_X)
\]
is bijective.
The tautological section therefore comes from a unique class, and $\mup$-gerbes up to equivalence are classified by their classes in $\HH^2(-,\mup)$.
This gives $\cG$.
The compatibility with $\iota$ holds because pullback of gerbes is functorial and the Kummer map commutes with pullback.

\textup{(iii)} Restricting the big fppf sheaf to $\Omega$-schemes gives the corresponding sheaf for $X_\Omega$, represented by $(\bH_X)_\Omega$ with identity component $\TT_\Omega\simeq\GGa$.
By \textup{(ii)}, the identification $\bH_X(\Omega)=\HH^2(X_\Omega,\mup)$ of \cite[(2.2.5.1)]{BLtw} sends $c$ to the class of $\cG_c$.
The surface $X_\Omega$ is supersingular by \Cref{prelim:base-change}, so $c\mapsto\alpha_c$ is the surjective homomorphism $U^2(X_\Omega,\mup)\to\Br(X_\Omega)$ of \cite[Theorem~2.2.6]{BLtw}, where $U^2(X_\Omega,\mup)=\TT(\Omega)$.

By the Kummer sequence and \Cref{prelim:base-change}, the kernel is the intersection of $\TT(\Omega)$ with the image of the finite subgroup $P\subset\bH_X(\bk)$ of Kummer classes of line bundles on $X$.
Since $\TT$ is open and closed in $\bH_X$, the intersection is $Z_X=\TT(\bk)\cap P$ for every $\Omega$.
This gives nonvanishing at transcendental parameters (compare \cite[proof of Theorem~6.2]{LieTw} for the geometric generic point).
Equivariance follows from \textup{(ii)}.
\end{proof}

For $c\in\Omega$ we put $X_c=X(\alpha_c)$.
By \Cref{prelim:twistor}\textup{(iii)} applied to $X=J$, the Artin--Tate twists of $J$ over $\Omega$ are the surfaces $X_c$ with $c\in\Omega\smallsetminus Z_J$.
The gerbe $\cG$ attaches a class to every point of the twistor line, not only to its $\bk$-points.
A parameter transcendental over $\bk$ is a geometric point of $\TT$, and the following case is the one used in the proof of \Cref{thm:main}.

\begin{corollary}\label{prelim:generic-twist}
Let $\Omega=\kab$, so that the coordinate $a$ defines the geometric generic point $\Spec\Omega\to\TT$.
Then $\alpha_a\neq0$, and the \emph{geometric generic Artin--Tate twist} $X_a$ over $\Omega$ is a supersingular $K3$ surface of Artin invariant $\sigma+1$.
In particular, $\sigma(X_a)=10$ if $\sigma=9$. 
The same holds for $X_c$ for every $c\in\Omega$ that is transcendental over $\bk$.
We call the generic fiber of $X_a$, an $E_{\Omega(t)}$-torsor, the \emph{geometric generic torsor}.
\end{corollary}
\begin{proof}
By \Cref{prelim:twistor}\textup{(iii)}, $\alpha_c\neq0$ for such $c$, and in particular for $c=a$.
The remaining assertions follow from \Cref{prelim:artin-tate}.
\end{proof}

\begin{corollary}\label{prelim:generic-descent}
Let $\Omega=\kab$.
If the geometric generic Artin--Tate twist $X_a$ is unirational over $\Omega$, then $J$ is unirational over $\bk$.
More generally, $J$ is unirational as soon as $J_\Omega$ is unirational for some algebraically closed extension $\Omega/\bk$.
\end{corollary}
\begin{proof}
By \Cref{prelim:generic-twist} and \Cref{prop:AT twistor family}\textup{(i)}, the surface $X_a$ dominates $J_\Omega$, so it suffices to prove the last assertion.
Let $\phi\colon\PP^2_\Omega\dashrightarrow J_\Omega$ be dominant.
We can spread it out over a finitely generated $\bk$-subalgebra $A\subset\Omega$: there are a dense open $U\subset\PP^2_A$ and a morphism $\Phi\colon U\to J\times_\bk\Spec A$ over $A$ whose generic fiber restricts $\phi$.
Since $\phi$ is dominant, so is $\Phi$, and by Chevalley's theorem its image contains a dense open subset $W$.
Its complement $Z$ has generic fiber of dimension at most one, so for $s$ in a dense open subset of $\Spec A$ the fiber $Z_s$ is a proper closed subset of $J$, and $U_s\neq\emptyset$.
For such a closed point $s$, the image of $U_s$ contains $W_s\neq\emptyset$.
Thus $\Phi_s\colon\PP^2_\bk\dashrightarrow J$ is dominant.
\end{proof}

\begin{remark}[The twistor family]\label{prelim:twistor-family}
We use the language of Bragg and Lieblich \cite{BLtw}.
The gerbe $\cG$ of \Cref{prelim:twistor} is for $X=J$ the tautological family of $\mup$-gerbes on $J$, the basic example of a twistor family \cite[Section~2.2]{BLtw}.
In the twisted period domain, $\TT$ corresponds to a twistor line in the sense of \cite[Definition~3.2.2]{BLtw}. This follows from \cite[Propositions~3.1.11 and~3.5.9]{BLtw}, as explained in the proof of Proposition~3.5.10 there.
For $c\in\Omega$ with $\alpha_c\neq0$, the Artin--Tate twist $X_c$ is the coarse moduli space of stable $\cG_c$-twisted sheaves of Mukai vector $(0,f,1)$, with a generic polarization as in the proof of \Cref{prelim:artin-tate}.
In this sense the surfaces $X_c$ form the Artin--Tate twistor family of $J$ over $\TT$.
\end{remark}

\section{Separable rational multisections}\label{sec:separable rational multisections}
\subsection{Splitting fields of Artin--Tate torsors}
We keep the standing notation of \Cref{sec:artin-tate}.
A torsor \emph{splits} over an extension of its base field if it has a point over that extension.
A multisection $M$ of a fibration over $\PP^1_\Omega$ is \emph{separable} or \emph{purely inseparable} if the extension $\Omega(M)/\Omega(t)$ of function fields has that property.
An elliptic curve over $k(t)$ is isotrivial if it becomes a curve defined over $k$ after a finite extension of $k(t)$; equivalently, its $j$-invariant function is constant.
Otherwise, we call it non-isotrivial.
\begin{lemma}\label{prelim:fixed-cover}
  Let $l$ be an algebraically closed field, let $E$ be an elliptic curve over $l(t)$, and let $L/l(t)$ be a finite separable extension.
  Then the group
  \[
    \ker\left(\HH^1(l(t),E)\longrightarrow\HH^1(L,E)\right)
  \]
  is finite.
  In particular, only finitely many classes in $\Sha(E)$ split over $L$.
\end{lemma}
\begin{proof}
  Put $K=l(t)$.
  Enlarging $L$ only enlarges the kernel, so we may replace $L$ by its Galois closure in $K^{\sep}$ and assume that $L/K$ is Galois, with group $G$.
  Inflation--restriction for the normal subgroup $\Gal(K^{\sep}/L)$ of $\Gal(K^{\sep}/K)$ acting on $E(K^{\sep})$ identifies the kernel with $\HH^1(G,E(L))$ \cite[Ch.~VII, \S6, Proposition~4]{SerreLF}.

  Since $l$ is algebraically closed, $L/l$ is a finitely generated regular extension.
  Let $\tau\colon B_L\to E_L$ be the $L/l$-trace of $E_L$, where $B$ is an abelian variety over $l$.
  The map $\tau$ is injective on $L$-points \cite[Theorem~6.12 and the beginning of \S7]{Conrad}.
  By the Lang--N\'eron theorem \cite{LN}, \cite[Theorem~7.1]{Conrad}, the quotient $Q=E(L)/A$ by the subgroup $A=\tau(B(l))\cong B(l)$ is finitely generated.
  Since $l$ is algebraically closed, multiplication by any integer $n\geq1$ is surjective on $B(l)$, with finite kernel $B[n](l)$.
  So $A$ is divisible, and $A[n]$ is finite for every $n\geq1$.
  The image in $Q$ of a divisible subgroup of $E(L)$ is a divisible subgroup of a finitely generated abelian group, hence zero.
  Thus $A$ is the largest divisible subgroup of $E(L)$.
  In particular, $A$ is $G$-stable, and \[0\longrightarrow A\longrightarrow E(L)\longrightarrow Q\longrightarrow 0 \] is an exact sequence of $G$-modules.

  Let $n=\lvert G \rvert$.
  The groups $\HH^1(G,A)$ and $\HH^1(G,Q)$ are killed by $n$, and $\HH^1(G,Q)$ is finite because $Q$ is finitely generated \cite[Ch.~VIII, \S2, Corollaries~1 and~2 to Proposition~4]{SerreLF}.
  Since $A$ is divisible, the sequence $0\to A[n]\to A\xrightarrow{n}A\to0$ is exact.
  By the sequence after taking cohomologies, the $n$-torsion of $\HH^1(G,A)$ is the image of the finite group $\HH^1(G,A[n])$.
  Since $\HH^1(G,A) = \HH^1(G,A)[n]$, $\HH^1(G,A)$ itself is finite.
  The exact sequence \[\HH^1(G,A)\to\HH^1(G,E(L))\to\HH^1(G,Q)\] now implies that $\HH^1(G,E(L))$ is finite.
\end{proof}
If $E$ is moreover ordinary, there is a similar finiteness statement for the torsors that split over a purely inseparable extension.
\begin{lemma}\label{prelim:frobenius-finite}
  Keep the standing notation of \Cref{sec:artin-tate}, and assume that the geometric generic fiber of $\pi$ is an ordinary elliptic curve.
  Let $\Omega\supseteq\bk$ be algebraically closed and $e\geq1$.
  Then only finitely many classes in $\Br(J_\Omega)\simeq\Sha(E_{\Omega(t)})$ split over $\Omega(t)^{1/p^e}$.
\end{lemma}
\begin{proof}
  The base change $J_\Omega$ satisfies the same hypotheses over $\Omega$ by \Cref{prelim:base-change}, and ordinarity does not change under extension of the base field.
  Replacing $\bk$ by $\Omega$, we may assume $\Omega=\bk$.
  Put $K=\bk(t)$.

  Choose a nonempty affine open $U\subset\PP^1_\bk$ over which $\pi$ is smooth with ordinary fibers.
  This is possible because the Hasse invariant of the generic fiber is nonzero, so it vanishes at only finitely many points.
  Let $\cE=J|_U$, an elliptic curve over $U$ with generic fiber $E$.
  Let $\alpha$ be a class that splits over $K^{1/p^e}$, with torsor $T$, and put $\cT=X(\alpha)|_U$.
  By \Cref{prelim:artin-tate}, $\cT$ is \'etale locally on $U$ isomorphic to $\cE$, so it is smooth and proper over $U$ with genus-one fibers.
  Its relative Picard scheme $\Pic^0_{\cT/U}$ is an elliptic curve over $U$ with generic fiber $\Pic^0_T\cong E$, hence isomorphic to $\cE$, since an abelian scheme over $U$ is the N\'eron model of its generic fiber.
  As $\cT\to U$ is smooth, it has sections \'etale locally, so $\cT$ is an $\cE$-torsor extending $T$.

  Let $\Frob^e\colon\cE\to\cE^{(p^e)}$ be the $e$-fold relative Frobenius over $U$, and let $G_e=\ker\Frob^e$.
  It is a finite flat infinitesimal group scheme of order $p^e$, and the fppf sequence \[ 0\to G_e\to\cE\to\cE^{(p^e)}\to0 \] is exact.
  The relative Frobenius $\cT\to\cT^{(p^e)}$ is equivariant for $\Frob^e$, so $\cT^{(p^e)}$ is the pushout of $\cT$ along $\Frob^e$.
  Identifying $K^{1/p^e}$ with $K$ by $y\mapsto y^{p^e}$, a $K^{1/p^e}$-point of $T$ becomes a $K$-point of $T^{(p^e)}$.
  Since $\cT^{(p^e)}$ is proper over the regular curve $U$, this point extends to a section over $U$.
  Hence $\Frob^e_*[\cT]=0$, and $[\cT]$ lies in the image of $\HH^1(U,G_e)\to\HH^1(U,\cE)$.
  Restricting to the generic point, $\alpha$ lies in the image of the composite $\HH^1(U,G_e)\to\HH^1(K,E)$, which does not depend on $\alpha$.
  It suffices to show that $\HH^1(U,G_e)$ is finite.

  Since the fibers of $\cE$ are ordinary, the Cartier dual of $G_e$ is the kernel of the $e$-fold Verschiebung, which is \'etale and cyclic of order $p^e$.
  Choose a connected finite \'etale cover $V\to U$ trivializing it.
  Then $G_e|_V\cong\mu_{p^e}$.
  The restriction $\HH^1(U,G_e)\to\HH^1(V,G_e)$ is injective.
  Indeed, $G_e$ is infinitesimal, so it has no nonzero sections over a reduced scheme.
  A $G_e$-torsor over $U$ that is trivial over $V$ therefore has a unique section over $V$.
  Its two pullbacks to the reduced scheme $V\times_UV$ agree, so it descends to a section over $U$.

  The fppf Kummer sequence gives an exact sequence
  \[
    0\longrightarrow\Oh(V)^\times/\Oh(V)^{\times p^e}\longrightarrow\HH^1(V,\mu_{p^e})\longrightarrow\Pic(V)[p^e]\longrightarrow0.
  \]
  The curve $V$ is smooth, connected and affine over $\bk$.
  Let $\bar V$ be its smooth compactification and $B\subset\Pic(\bar V)$ the finitely generated subgroup generated by the points of $\bar V\smallsetminus V$.
  The group $\Oh(V)^\times/\bk^\times$ embeds into the free abelian group on $\bar V\smallsetminus V$, and $\bk^\times$ is divisible, so the left-hand group is finite.
  Applying the snake lemma to multiplication by $p^e$ on $0\to B\to\Pic(\bar V)\to\Pic(V)\to0$ gives an exact sequence
  \[
    \Pic(\bar V)[p^e]\longrightarrow\Pic(V)[p^e]\longrightarrow B/p^eB.
  \]
  Here $\Pic(\bar V)[p^e]$ is the $p^e$-torsion of the Jacobian of $\bar V$, which is finite, and $B/p^eB$ is finite.
  Hence $\Pic(V)[p^e]$ is finite, and so are $\HH^1(V,\mu_{p^e})$ and $\HH^1(U,G_e)$.
\end{proof}
\begin{remark}
The statement of \Cref{prelim:frobenius-finite} fails when the generic elliptic curve is geometrically supersingular.
In fact, one can show that the map
\[ \Sha(E_{\Omega(t)}) \longrightarrow \HH^1\bigl(\Omega(t)^{1/p^e},E\bigr)\]
is then zero for every $e \geq 2$, because multiplication by $p$ on a supersingular elliptic curve factors through the $2$-fold Frobenius.
Thus infinitely many torsors split over $\Omega(t)^{1/p^e}$ for $e \geq 2$.
The property used in the ordinary case is that the kernel $G_e$ is a connected finite group scheme of multiplicative type.
\end{remark}
\begin{proposition}\label{prelim:no-inseparable}
  Keep the standing notation of \Cref{sec:artin-tate}, and assume that the geometric generic fiber of $\pi$ is an ordinary elliptic curve.
  Let $\Omega\supseteq\bk$ be an algebraically closed extension.
  \begin{enumerate}[label=\textup{(\roman*)},leftmargin=2.2em]
    \item For every $c\in\Omega\smallsetminus\bk$, the Artin--Tate twist $X_c\to\PP^1_\Omega$ has no purely inseparable multisection.
    \item For every $e\geq0$, only finitely many $c\in\bk$ have the property that $X_c$ has a purely inseparable multisection of degree $p^e$.
  \end{enumerate}
\end{proposition}
\begin{proof}
  For $e\geq0$, let \[S_e(\Omega)\subset\Omega\] be the set of $c$ such that $\alpha_c$ splits over $\Omega(t)^{1/p^e}$.
  The map $c\mapsto\alpha_c$ has the finite kernel $Z_J$ by \Cref{prelim:twistor}\textup{(iii)}.
  The set $S_0(\Omega)=Z_J$ is finite, since a torsor with an $\Omega(t)$-point is trivial.
  For $e\geq1$, \Cref{prelim:frobenius-finite} shows that $S_e(\Omega)$ is the preimage of a finite set, hence finite.

  Let $M\subset X_c$ be a purely inseparable multisection.
  Its degree is $[\Omega(M):\Omega(t)]=p^e$ for some integer $e\geq 0$, and every element of $\Omega(M)$ has a $p^e$th power in $\Omega(t)$.
  Thus $\Omega(M)\subseteq\Omega(t)^{1/p^e}$.
  Since $\Omega$ is perfect, $\Omega(t)^{1/p^e}=\Omega(t^{1/p^e})$ has degree $p^e$ over $\Omega(t)$, so $\Omega(M)=\Omega(t)^{1/p^e}$.
  The generic point of $M$ is a closed point of the generic fiber of $X_c$ with residue field $\Omega(M)$.
  It is a point of the torsor with class $\alpha_c$ over $\Omega(t)^{1/p^e}$, so $c\in S_e(\Omega)$.
  Applying it for $\Omega=\bk$, we can conclude \textup{(ii)}.

  For \textup{(i)}, let $\tau\in\Aut(\Omega/\bk)$, extended to $\Omega(t)$ by $\tau(t)=t$.
  It maps $\Omega(t)^{1/p^e}$ onto itself.
  By the equivariance in \Cref{prelim:twistor}\textup{(iii)} and the canonicity of $\Br(J_\Omega)\simeq\Sha(E_{\Omega(t)})$, base change along $\tau$ carries the torsor with class $\alpha_c$ to the torsor with class $\alpha_{\tau(c)}$.
  Hence $S_e(\Omega)$ is stable under $\Aut(\Omega/\bk)$.
  An element $c\in\Omega\smallsetminus\bk$ is transcendental over $\bk$, since $\bk$ is algebraically closed.
  For $\lambda\in\bk$, the automorphism $c\mapsto c+\lambda$ of $\bk(c)$ extends to $\overline{\bk(c)}$.
  Acting trivially on a transcendence basis of $\Omega$ over $\overline{\bk(c)}$, it extends further to $\Omega$.
  So the orbit of $c$ contains $c+\bk$, which is infinite.
  A finite stable set contains no point with infinite orbit, so $S_e(\Omega)\subset\bk$ for every $e$.
  Hence $c\notin S_e(\Omega)$ for every $e$, and $X_c$ has no purely inseparable multisection.
\end{proof}

A non-isotrivial elliptic curve $E$ is automatically ordinary: its $j$-invariant is transcendental over $\bk$, whereas supersingular $j$-invariants lie in $\FF_{p^2}$.
The Artin invariant of $J$ plays no role, and $p=3$ is allowed.
By \textup{(ii)}, the parameters $c\in\bk$ for which $X_c$ has a purely inseparable multisection form a countable set, so a very general $\bk$-point of the twistor line has none when $\bk$ is uncountable.

\begin{remark}\label{prelim:BL-comparison}
  Suppose $p\geq5$ and that for each $e\geq1$ the Frobenius pullback map $\gamma_e(\pi)$ from $\HH^2(J,\Oh_J)$ to the corresponding cohomology of the resolved Frobenius pullback is injective \cite[(4.5.1)]{BLpoints}. This holds, for instance, when $\pi$ is $\infty$-Frobenius split \cite[Lemma~4.6]{BLpoints}.
  Completing $\TT$ at the origin gives the universal formal Brauer class by \cite[Lemma~2.2.3]{BLtw}.
  The geometric generic member over $\overline{\bk((a))}$ then has no purely inseparable multisection by \cite[Corollary~6.7 and proof of Corollary~6.9]{BLpoints}.
  A purely inseparable multisection of $X_a$ over $\overline{\bk(a)}$ would persist after extending scalars to $\overline{\bk((a))}$ along an embedding extending $\bk(a)\hookrightarrow\bk((a))$.
  Injectivity of all $\gamma_e(\pi)$ forces an ordinary generic fiber by the proof of \cite[Proposition~6.10]{BLpoints}. For ordinary generic fibers, \Cref{prelim:no-inseparable} gives the conclusion in every odd characteristic.
  If $p\geq5$ and the geometric generic fiber is supersingular, then $\pi$ is isotrivial and every twist has a purely inseparable multisection by \cite[Proposition~6.10(i)]{BLpoints}.
\end{remark}

The separable analogue of \Cref{prelim:no-inseparable}\textup{(i)} concerns covers of the base that are defined over $\bk$: such a cover splits only the twists over $\bk$-points of the twistor line.

\begin{corollary}\label{prelim:moving-covers}
  Keep the standing notation of \Cref{sec:artin-tate}.
  Let $\Omega\supseteq\bk$ be algebraically closed, and let $C\to\PP^1_\bk$ be a finite separable cover by an integral curve.
  Then the set of $c\in\Omega$ such that $\alpha_c$ splits over the function field of $C_\Omega$ is finite and contained in $\bk$.
  In particular, $C_\Omega\to\PP^1_\Omega$ does not split the geometric generic torsor.
\end{corollary}
\begin{proof}
  Since $\bk$ is algebraically closed, $C_\Omega$ is integral, and its function field $L=\Omega(t)\otimes_{\bk(t)}\bk(C)$ is a finite separable extension of $\Omega(t)$.
  Let $S\subset\Omega$ be the set in question.
  By \Cref{prelim:fixed-cover}, applied with $l=\Omega$ to $E_{\Omega(t)}$, only finitely many classes in $\Br(J_\Omega)\simeq\Sha(E_{\Omega(t)})$ split over $L$.
  The map $c\mapsto\alpha_c$ has the finite kernel $Z_J$ by \Cref{prelim:twistor}\textup{(iii)}, so $S$ is finite.
  Let $\tau\in\Aut(\Omega/\bk)$, extended to $\Omega(t)$ by $\tau(t)=t$ and to $L$ by acting trivially on $\bk(C)$.
  As in the proof of \Cref{prelim:no-inseparable}\textup{(i)}, base change along $\tau$ carries the torsor with class $\alpha_c$ to the torsor with class $\alpha_{\tau(c)}$.
  Hence $S$ is stable under $\Aut(\Omega/\bk)$.
  Every $c\in\Omega\smallsetminus\bk$ has an infinite orbit under $\Aut(\Omega/k)$, so $S \subset\bk$ by the similar argument again.
  For $\Omega=\kab$, the parameter $a$ is not in $\bk$, which proves the last assertion.
\end{proof}

\subsection{Descent to the twistor line}
\label{sec:roots-descent}

In this subsection $J$ is as in the standing notation of \Cref{sec:artin-tate}.
We put $L=\bk(a)$, the function field of the twistor line $\TT$, and $\Omega=\overline L$.

Let $\cG_a\to J_L$ be the pullback of the universal gerbe $\cG$ along the generic point $a\colon\Spec L\to\TT$.
Its base change to $\Omega$ is the gerbe of the geometric generic twist by \Cref{prelim:twistor}\textup{(ii)}.
By \Cref{prelim:twistor-family}, $X_a\to\PP^1_\Omega$ is the moduli space of stable twisted sheaves of Mukai vector $(0,f,1)$. The polarization comes from $J$ by \Cref{prelim:base-change}, and the fiber class is defined over $L$.
We use the relative moduli problem over $L$ of \cite[Definition~4.1.17]{BLtw}, which tests stability and the Mukai vector on geometric fibers.
Its base change to $\Omega$ agrees with the fiberwise moduli problem with numerical parameter $1$ by \cite[Definition~4.3.2(2) and Lemma~4.3.5]{BLtw}.
The stack of stable twisted sheaves is algebraic, as can be seen by representing the Brauer class of $\cG_a$ by an Azumaya algebra on $J_L$ and using the stack of coherent modules over that algebra.
Rigidifying the scalar automorphisms gives an algebraic stack with trivial inertia \cite[\href{https://stacks.math.columbia.edu/tag/04V2}{Tag~04V2}]{stacks-project}, hence an algebraic space $X\to\PP^1_L$ \cite[\href{https://stacks.math.columbia.edu/tag/04SZ}{Tag~04SZ}]{stacks-project}. The rigidification is an fppf gerbe, so $X$ represents the fppf sheafification of isomorphism classes.
This moduli space commutes with extension of the base field.
Since $X_\Omega\cong X_a$ is a $K3$ surface, $X$ is smooth and proper over $L$.
Properness makes $X$ decent and of finite type, and every finite set of points of $X_\Omega$ lies in an affine open. Thus $X$ is a scheme by \cite[\href{https://stacks.math.columbia.edu/tag/0B88}{Tag~0B88}]{stacks-project}.
It is quasi-projective because $X_\Omega$ is \cite[\href{https://stacks.math.columbia.edu/tag/0BDE}{Tag~0BDE}]{stacks-project}, hence projective by properness.

Consider the following condition for the fibers of $\pi \colon J \to \PP^1_k$:
\begin{equation}\label{eq:star}
  \text{All fibers of $\pi$ are reduced and their component groups have order prime to $p$.} \tag{{$\star$}} 
\end{equation}
Under \eqref{eq:star}, the geometric Picard group has trivial Galois action, and the smooth rational curves descend to $L$.

\begin{lemma}\label{roots:descent}
Assume that all fibers of $\pi \colon J \to \PP^1_k$ are reduced and that their component groups have order prime to $p$.
Let $X\to\PP^1_L$ be the model of its geometric generic Artin--Tate twist over $L=k(a)$, and put $\Omega=\overline L$.
Then $\Gal(L^{\sep}/L)$ acts trivially on $\Pic(X_\Omega)$, and every smooth rational curve on $X_\Omega$ descends to a curve on $X$ isomorphic to $\PP^1_L$.
\end{lemma}
\begin{proof}
Put $G=\Gal(L^{\sep}/L)$, which acts on $\Omega$ by unique extension.
The map $\Pic(X_{L^{\sep}})\to\Pic(X_\Omega)$ is an isomorphism \cite[Lemma~3.1]{BLvL}, so $G$ acts on $\Pic(X_\Omega)$ through a finite quotient.

Let $T$ be the generic fiber of $X$, a torsor under $E_{L(t)}$ parametrizing invertible twisted sheaves of a fixed degree on the restriction of $\cG_a$ to $E_{L(t)}$.
For such a sheaf $\mathcal L$, the power $\mathcal L^{\otimes p}$ has trivial inertia action and comes from an invertible sheaf of fixed degree $d$ on $E_{L(t)}$.
The class of $\mathcal L^{\otimes p}(-dO)$, where $O$ is the origin, gives a $[p]$-map $m\colon T\to E_{L(t)}$.
Indeed the construction is compatible with base change and tensoring $\mathcal L$ by a degree-zero line bundle changes the class by $p$ times that line bundle.
Over $\PP^1_\Omega$, the relative smooth locus of $J_\Omega$ is the N\'eron model of its generic fiber.
Its mapping property extends $m$ to a morphism from the relative smooth locus of $X_\Omega$.
The extension is $G$-equivariant by uniqueness.
After an \'etale-local identification of the two fibrations, the extension acts on the component groups as multiplication by $p$ followed by a translation, so it is bijective on components by \eqref{eq:star}.
Every component meets the relative smooth locus since the fibers are reduced.
Thus $m$ matches the components of $X_\Omega$ equivariantly with those of $J_\Omega$, which are defined over $\bk$, and $G$ fixes the subgroup $V\subset\Pic(X_\Omega)$ generated by vertical classes.

Restriction to the generic fiber identifies $\Pic(X_\Omega)/V$ with $\Pic(T_{\Omega(t)})$, since divisors on the generic fiber extend by closure on the smooth surface $X_\Omega$.
Its degree-zero subgroup is $E(\Omega(t))$, since $\Br(\Omega(t))=0$ by Tsen's theorem.
This equals $E(\bk(t))$, since both are the quotients of $\Pic(J_\Omega)=\Pic(J)$ (\Cref{prelim:base-change}) by the vertical classes and the zero section.
Thus $G$ acts trivially on $V$, on the degree-zero subgroup of $\Pic(T_{\Omega(t)})$, and on the degree quotient.
For $\sigma\in G$ this gives $(\sigma-1)^3=0$ on $\Pic(X_\Omega)$, while $\sigma$ has finite order on this group.
Hence $\sigma$ acts as the identity on $\Pic(X_\Omega)\otimes\QQ$, and therefore on the torsion-free group $\Pic(X_\Omega)$ (see \cite[Corollaries~5.21--5.22]{LieTw}).

Let $C\subset X_\Omega$ be a smooth rational curve.
By \cite[Corollary~3.2]{BLvL}, $C$ descends to a curve $C'\subset X_{L^{\sep}}$.
Since $C'^2=-2$, it is the unique effective divisor in its class.
Its class is $G$-fixed, so $C'$ descends to a smooth conic $C_0\subset X$.
It has an $L$-point because $L$ is $C_1$ by Tsen's theorem, so $C_0\cong\PP^1_L$.
\end{proof}

When $p\geq5$ and $\sigma(J)=9$, the condition \eqref{eq:star} holds for every Jacobian elliptic fibration on $J$, by \Cref{lem:Kodaira}.
In characteristic three a fiber of type $\KI_3$ or $\KIV$ has component group $\ZZ/3$, and \eqref{eq:star} is checked directly for the fibration of \Cref{pencils:three}.

\begin{proposition}\label{prop:existence of seperable rational multisection}
Let $\pi\colon J\to\PP^1_\bk$ be a supersingular Jacobian elliptic $K3$ surface with $\sigma(J)=9$ in odd characteristic $p$, whose geometric generic fiber is ordinary.
\begin{enumerate}[label=\textup{(\roman*)},leftmargin=2.2em]
  \item The geometric generic Artin--Tate twist $X_a\to\PP^1_\Omega$ carries a smooth rational multisection of degree $p$, and every multisection of degree $p$ of $X_a$ is separable over $\PP^1_\Omega$.
  \item If $\pi$ satisfies \eqref{eq:star}, let $X\to\PP^1_{\bk(a)}$ be the model of $X_a$ of \Cref{roots:descent}.
    Then every smooth rational multisection of degree $p$ of $X_a$ descends to a separable multisection of $X\to\PP^1_{\bk(a)}$ of degree $p$, isomorphic to $\PP^1_{\bk(a)}$.
\end{enumerate}
\end{proposition}
\begin{proof}
(i) By \Cref{prelim:generic-twist}, $X_a$ is a supersingular $K3$ surface of Artin invariant ten.
By \Cref{prelim:artin-tate}, $X_a\to\PP^1_\Omega$ is \'etale locally isomorphic to $\pi_\Omega$, so it is an elliptic fibration.
Hence \Cref{roots:reflection}\textup{(ii)} and \textup{(iii)}, applied over the algebraically closed field $\Omega$, give a smooth rational multisection of degree $p$.
Let $M$ be any multisection of degree $p$ of $X_a$.
The separable closure of $\Omega(t)$ in $\Omega(M)$ has degree dividing the prime $p$, so $M$ is either separable or purely inseparable over $\PP^1_\Omega$.
Since $a\notin\bk$ and the geometric generic fiber of $\pi$ is ordinary, \Cref{prelim:no-inseparable}\textup{(i)} excludes the second case.

(ii) Let $M\subset X_a$ be a smooth rational multisection of degree $p$.
By \Cref{roots:descent}, it descends to a curve $M_0\subset X$ with $M_0\cong\PP^1_L$ and $M_{0,\Omega}=M$, where $L=\bk(a)$.
The morphism $M_0\to\PP^1_L$ is finite flat of degree $p$ by faithfully flat descent from $M\to\PP^1_\Omega$.
Its generic fiber is \'etale since its base change to $\Omega(t)$ is \'etale by \textup{(i)}.
\end{proof}

\section{Rational multisections and unirationality}\label{sec:rational multiplication and unirationality}

\subsection{Elliptic curves without separable \texorpdfstring{$p$}{p}-torsion}
Let $K$ be a field of characteristic $p$, not necessarily perfect (e.g. $k(t)$), and fix a separable closure $K^{\sep}$ of $K$.
Let $E$ be an elliptic curve over $K$.

If $E$ is supersingular, then its $p$-torsion subgroup scheme is connected, so $E(K^{\sep})[p]=0$.
Over a non-perfect field $K$ this can also happen when $E$ is ordinary, and it is controlled by the $j$-invariant of $E$.
\begin{lemma}\label{lem:no-p-torsion vs j-invariant}
If $E$ is ordinary, then 
\begin{equation}\label{eq:no separable p-torsion}
E(K^{\sep})[p] = 0 \quad \iff \quad j(E) \notin K^{p}\,.
\end{equation}
\end{lemma}
\begin{proof}
  
Suppose that $E(K^{\sep})[p]\neq0$.
A nonzero point generates an \'etale subgroup of order $p$, defined over a finite separable extension $L/K$.
The dual of the quotient isogeny is radicial of degree $p$ onto $E$, so it is the relative Frobenius of the quotient, and $E_L$ is a relative Frobenius twist.
Thus $j(E)$ is a $p$th power in $L$.
An element of $K$ that is a $p$th power in a separable extension is a $p$th power in $K$, since its $p$th root is both separable and purely inseparable over $K$; see \cite[Lemma 9.3]{LS10}.
Hence $j(E)\in K^p$.

For the converse, suppose that $j(E)\in K^p$.
By \cite[Proposition 3.3]{LS10} there is an elliptic curve $E'$ over $K$ whose Frobenius twist $E'^{(p)}$ is isomorphic to $E$ over $K$.
Since $E$ is ordinary, so is $E'$, and thus the Verschiebung $E'^{(p)} \to E'$ is an \'etale isogeny of degree $p$, whose kernel is a nontrivial subgroup of $E'^{(p)}(K^{\sep})[p]$.
Hence $E(K^{\sep})[p]\cong E'^{(p)}(K^{\sep})[p]\neq0$.
\end{proof}
Next, we show that the condition $E(K^{\sep})[p]=0$ is stable under extension of the field of constants of a function field.
\begin{lemma}\label{sweep:no-torsion}
 Let $E$ be an elliptic curve over $\bk(t)$, and let $\Omega\supseteq\bk$ be algebraically closed.
  The following conditions are equivalent:
  \begin{enumerate}
    \item $E(\Omega(t)^{\sep})[p]=0$.
    \item $E(\bk(t)^{\sep})[p]=0$.
    \item $E$ is geometrically supersingular or $j(E) \notin k(t)^{p}$.
  \end{enumerate}
\end{lemma}
\begin{proof}

 Geometric supersingularity does not change under extension of the constants.
For a perfect field $\bk'$, the $p$th powers in $\bk'(t)$ form the subfield $\bk'(t^p)$, which is the kernel of $d/dt$.
The $t$-derivative of $j(E)$ is the same in $\bk(t)$ and in $\Omega(t)$, so $j(E)\in\Omega(t)^p$ if and only if $j(E)\in\bk(t)^p$.
If $E$ is geometrically supersingular, then \textup{(1)} and \textup{(2)} hold.
Otherwise, applying \Cref{lem:no-p-torsion vs j-invariant} over $\bk(t)$ and over $\Omega(t)$ gives the equivalences.
\end{proof}

\subsection{Sweeping rational multisections along the twistor line}
We first record how to spread out a parametrization of a curve over $\bk(a)$.
\begin{lemma}\label{sweep:parameter}
Let $Y$ be a proper $\bk$-scheme and $\TT=\Spec\bk[a]$.
Every $\bk$-morphism $\psi\colon\Spec\bk(a,s)\to Y$ extends to a morphism $U\times\PP^1_\bk\to Y$ for some nonempty open $U\subset\TT$, where $s$ is an affine coordinate on $\PP^1_\bk$.
\end{lemma}
\begin{proof}
By properness of $Y$, the map extends over every local ring of the regular curve $\PP^1_{\bk(a)}$.
At each point of this generic fiber, its local ring equals the local ring of $\TT\times\PP^1_\bk$ at the same point.
Since $Y$ is of finite type over $\bk$, each local extension is defined on an open neighborhood of that point in $\TT\times\PP^1_\bk$.
The extensions agree on overlaps by separatedness and give a morphism on an open subset containing the generic fiber.
The projection $\TT\times\PP^1_\bk\to\TT$ is proper, so the image of the complement is closed and avoids the generic point of $\TT$.
Its complement is the required open subset $U$.
\end{proof}

Let $E$ be an elliptic curve over a field $K$ of characteristic $p$, and let $T$ be an $E$-torsor with a fixed effective divisor $D\subset T$ of degree $p$.
We do not assume that $D$ is \'etale over $K$.
For a field extension $L$ of $K$ and $x,y\in T(L)$, write $x-y\in E(L)$ for the difference in the torsor; it is given by a $K$-morphism $T\times_K T\to E$.
For a $K$-scheme $S$ and $x\in T(S)$, write $\Gamma_x\subset T_S$ for the image of the section $x$, an effective relative Cartier divisor of degree one, and $D_S\subset T_S$ for the pullback of $D$.

We identify $E$ with the Jacobian $\Pic^0_{T/K}$ of the smooth projective genus-one curve $T$, so that $x-y$ is the class of the divisor $[x]-[y]$ on $T_L$ \cite[Chapter~11, \S5.1]{Huy}; see \cite[\S2, Lemma~4 and Proposition~1]{Taelman} for the relative version over any base scheme.
The inverse $\Pic^0_{T/K}\to E$ is the homomorphism that \cite[Proposition~6.4]{MilneJac} attaches to the difference map $T\times_K T\to E$, which vanishes on the diagonal.
Explicitly, the isomorphism $E\times_K T\to T\times_K T$, $(e,y)\mapsto(y+e,y)$, turns $(x,y)\mapsto\Oh_{T_S}(\Gamma_x-\Gamma_y)$ into a morphism $\beta\colon E\times_K T\to\Pic^0_{T/K}$.
Over $\bar K$, choose $o\in T(\bar K)$ and identify $T_{\bar K}$ with $E_{\bar K}$ by $x\mapsto x-o$.
Abel's theorem on the elliptic curve $E_{\bar K}$ gives $\beta(e,y)=\Oh([o+e]-[o])$, which does not depend on $y$, and $e\mapsto\Oh([o+e]-[o])$ is an isomorphism of group schemes $E_{\bar K}\to\Pic^0_{T_{\bar K}}$.
Hence the two pullbacks of $\beta$ to $E\times_K T\times_K T$ agree, since they agree on $\bar K$-points and $E\times_K T\times_K T$ is geometrically reduced.
By fpqc descent along $T\to\Spec K$, the morphism $\beta$ factors through an isomorphism of group schemes $E\to\Pic^0_{T/K}$.
\begin{definition}[Abel--Jacobi map]\label{defn:Abel-Jacobi map}
  The \emph{Abel--Jacobi map} $m_D\colon T\to E$ is the $K$-morphism that sends $x\in T(S)$, for a $K$-scheme $S$, to the class of $\Oh_{T_S}(p\Gamma_x-D_S)$ in $\Pic^0_{T/K}(S)=E(S)$.
\end{definition}
The invertible sheaf $\Oh_{T_S}(p\Gamma_x-D_S)$ has degree zero on every fiber of $T_S\to S$, and its formation commutes with base change in $S$, so $m_D$ is a well-defined morphism by the Yoneda lemma.
Let $L\supseteq K$ be a field over which $D_L=[d_1]+\dots+[d_p]$ with $d_1,\dots,d_p\in T(L)$, not necessarily distinct; for example, $L$ algebraically closed.
Since $p\Gamma_x-D_S=\sum_{i=1}^p(\Gamma_x-\Gamma_{d_i})$ for every $L$-scheme $S$ and $x\in T(S)$, the base change of $m_D$ to $L$ is the morphism
\[
  T_L\longrightarrow E_L,\qquad x\longmapsto\sum_{i=1}^{p}(x-d_i).
\]
When $D$ is \'etale over $K$, one may take for $L$ a finite Galois extension $K'\subset K^{\sep}$ of $K$, and $m_D$ is the Galois descent of this morphism, as $\Gal(K'/K)$ permutes the $d_i$.
\begin{lemma}\label{sweep:abel}
Let $D \subset T$ be an effective divisor of degree $p$.
\begin{enumerate}[label=\textup{(\roman*)},leftmargin=2.2em]
\item The formation of $m_D$ commutes with extension of $K$, and $m_D(x)-m_D(y)=[p](x-y)$ for every field $L\supseteq K$ and all $x,y\in T(L)$.
\item If $D$ is \'etale over $K$ and $E(K^{\sep})[p]=0$, then $m_D$ restricts to a closed immersion $D\to E$.
\end{enumerate}
\end{lemma}
\begin{proof}
(i) The definition is compatible with extension of $K$.
The difference $m_D(x)-m_D(y)$ is the class of $(p[x]-D_L)-(p[y]-D_L)=p([x]-[y])$, which is $[p](x-y)$.

(ii) Since $D$ is \'etale of degree $p$, there is a finite Galois extension $K'\subset K^{\sep}$ of $K$ over which $D$ splits into $p$ distinct points $d_1,\dots,d_p\in T(K')$.
For $i\neq j$, \textup{(i)} gives $m_D(d_i)-m_D(d_j)=[p](d_i-d_j)\neq0$, since $E(K^{\sep})[p]=0$.
So $m_D$ maps the $p$ reduced points of $D_{K'}$ injectively to rational points of $E_{K'}$, which is a closed immersion.
Closed immersions descend along field extensions \cite[\href{https://stacks.math.columbia.edu/tag/02L6}{Tag~02L6}]{stacks-project}.
\end{proof}
The \'etaleness in \textup{(ii)} cannot be dropped.
If $D=p[d]$ for some $d\in T(K)$, then $m_D(x)=[p](x-d)$, and $x\mapsto x-d$ maps $D$ isomorphically onto the subscheme $\Spec\Oh_{E,0}/\mathfrak{m}_0^p$ of $E$.
The multiplication $[p]$ factors through the relative Frobenius of $E$, which is constant on this subscheme, so $m_D$ contracts the non-reduced scheme $D$ to the point $0$.

The Abel--Jacobi map transfers multisections of an Artin--Tate twist into curves on $J$.
Let $X\to\PP^1_\Omega$ be an Artin--Tate twist of $J$ over an algebraically closed field $\Omega\supseteq\bk$, and let $M\subset X$ be a separable multisection of degree $p$.
Its generic point $D$ is a closed point of degree $p$ of the generic fiber of $X$, an $E_{\Omega(t)}$-torsor, and $D$ is \'etale over $\Omega(t)$.
Let $m_D(M)$ denote the Zariski closure of the point $m_D(D)$ in $J_{\Omega}$.
We call it the \emph{Abel--Jacobi image} of $M$.

In the \cref{fig:twistor-sweep} of the introduction, the Abel--Jacobi images of multisections of the twists $X_c$ sweep out $J$, as $c$ moves along the twistor line.
At the generic point of the twistor line, this means that the Abel--Jacobi image is a curve on $J_\Omega$ that does not come from a curve on $J$.

\begin{proposition}[Movable Abel--Jacobi images]\label{sweep:abel-image}
Let $J\to\PP^1_\bk$ be a supersingular Jacobian elliptic $K3$ surface in odd characteristic, whose generic elliptic curve $E/\bk(t)$ satisfies $E(\bk(t)^{\sep})[p]=0$.
Let $\Omega=\kab$, and let $M$ be a separable multisection of degree $p$ of the geometric generic Artin--Tate twist $X_a\to\PP^1_\Omega$, with generic point $D$.
\begin{enumerate}[label=\textup{(\roman*)},leftmargin=2.2em]
\item The Abel--Jacobi image $m_D(M)$ is a multisection of $J_\Omega\to\PP^1_\Omega$ of degree $p$, and $m_D$ maps $D$ isomorphically onto its generic point.
  In particular, $M$ and $m_D(M)$ are birational over $\PP^1_\Omega$.
\item The Abel--Jacobi image is not defined over $\bk$, i.e., there is no curve $C\subset J$ with $m_D(M)=C_\Omega$.
\end{enumerate}
\end{proposition}
By \textup{(i)}, the two multisections have isomorphic function fields over $\Omega(t)$.
This field splits the torsor of $X_a$, whose class moves with $a$.
A cover defined over $\bk$ splits only the twists $X_c$ with $c\in\bk$ (\Cref{prelim:moving-covers}).
This explains informally why the Abel--Jacobi image moves.
\begin{proof}
By \Cref{sweep:no-torsion}, the assumption $E(k(t)^{\sep})[p]=0$ implies that $E(\Omega(t)^{\sep})[p]=0$.

(i) By \Cref{sweep:abel}\textup{(ii)}, $m_D$ maps $D$ isomorphically onto the closed point $m_D(D)$ of $E_{\Omega(t)}$.
Hence $m_D(M)$ is an integral curve with function field $\kappa(m_D(D))\cong\Omega(M)$ over $\Omega(t)$, and it dominates $\PP^1_\Omega$.
It contains no fiber, so it is quasi-finite and proper, hence finite, over $\PP^1_\Omega$.
It is flat because it is integral and dominates the regular curve $\PP^1_\Omega$.
Its degree is $[\Omega(M):\Omega(t)]=p$.

(ii) Suppose that $m_D(M)=C_\Omega$ for a curve $C\subset J$ (defined over $k$).
Then $C$ is integral and finite over $\PP^1_\bk$, since it is proper and dominates the base.
The generic fiber of $C_\Omega\to\PP^1_\Omega$ has coordinate ring $\Omega(t)\otimes_{\bk(t)}\bk(C)$, which by \textup{(i)} is the separable field extension $\Omega(M)/\Omega(t)$.
By faithfully flat descent of \'etaleness, $\bk(C)/\bk(t)$ is separable.
The generic fiber of $X_a$ has the point $D$ with residue field $\Omega(M)$, so the cover $C_\Omega\to\PP^1_\Omega$ splits the geometric generic torsor.
This contradicts \Cref{prelim:moving-covers}.
\end{proof}
\begin{theorem}\label{sweep:criterion}
Suppose that $p \geq 3$.
Let $\pi \colon J \to \PP_k^1$ be a Jacobian elliptic fibration on a supersingular $K3$ surface that satisfies \eqref{eq:star}, and let $E/k(t)$ be its generic fiber.
Assume that the following conditions hold.
\begin{enumerate}
   \item $E(k(t)^{\sep})[p] = 0$.
   \item The geometric generic Artin--Tate twist $X_a \to \PP^1_{\kab}$ has a \emph{separable} smooth rational multisection of degree $p$.
\end{enumerate}
Then there is a dominant rational map $\TT\times\PP^1_\bk\dashrightarrow J$ whose restriction to the geometric generic fiber of the first projection parametrizes the Abel--Jacobi image of the multisection in \textup{(2)}. In particular, $J$ is unirational.
\end{theorem}
\begin{proof}
Let $a$ be the coordinate on the twistor line $\TT$, and put $L=\bk(a)$, $K=L(t)$, and $\Omega=\overline L$.
Let $X\to\PP^1_L$ be the model of $X_a$ from \Cref{sec:roots-descent}, with generic fiber $T/K$.
By \Cref{roots:descent}, the multisection $M$ in hypothesis~(2) descends to a curve $M_0\subset X$ isomorphic to $\PP^1_L$.
Choose a coordinate $s$ on $M_0$ over $L$.
The generic fiber $D_0=M_0\times_{\PP^1_L}\Spec K$ of $M_0\to\PP^1_L$ is a closed point of $T$ with residue field $L(s)$.
It is finite \'etale of degree $p$ over $K$ by descent from the generic fiber $D$ of $M$ over $\Omega(t)$.

The Abel--Jacobi map $m_{D_0}\colon T\to E_K$ is defined over $K$ by \Cref{defn:Abel-Jacobi map}.
Its restriction to $D_0$ gives the composite
\[
   \psi\colon \Spec\bk(a,s)=D_0\xrightarrow{\;m_{D_0}\;}E_K\longrightarrow J.
\]
By \Cref{sweep:parameter}, $\psi$ extends to a morphism on $U\times\PP^1_\bk$ for some nonempty open $U\subset\TT$. It gives the rational map
\begin{equation}\label{eq:sweep-map}
    \Phi\colon\TT\times\PP^1_\bk\dashrightarrow J
\end{equation}
whose restriction to the generic slice $\PP^1_L\cong M_0$ is induced by $m_{D_0}$.
The composite with $\pi$ is given by $t\in L(s)$, which is nonconstant since $M_0$ is a multisection.

After extending scalars to $\Omega(t)$, the point $D_0$ becomes $D$, and $m_{D_0}$ becomes $m_D$ by \Cref{sweep:abel}\textup{(i)}.
Thus the image of the geometric generic slice is the Abel--Jacobi image $m_D(M)$.
If $\Phi$ were nondominant, the closure of its image would be an integral curve $C\subset J$ over $\bk$.
Then $\psi$ would factor through $C$, and extending the slice to $\Omega$ would give $m_D(M)\subset C_\Omega$.
Since $\bk$ is algebraically closed, $C_\Omega$ is integral, so the two curves would coincide.
This contradicts \Cref{sweep:abel-image}\textup{(ii)} and proves dominance.
\end{proof}

\begin{corollary}\label{cor:criterion p>5}
 Suppose that $p \geq 5$. Let $X$ be a supersingular $K3$ surface with $\sigma(X) =9$.
Assume that $X$ admits a Jacobian elliptic fibration whose generic elliptic curve has $j$-invariant function $j(E) \notin k(t)^p$. Then $X$ is unirational.
\end{corollary}
\begin{proof}
Let $\pi\colon X=J\to\PP^1_\bk$ be such a fibration, with generic elliptic curve $E/\bk(t)$, and let $\Omega=\kab$.
As $\sigma(J)=9$ and $p\geq5$, the fibration $\pi$ satisfies \eqref{eq:star}, as explained before \Cref{prop:existence of seperable rational multisection}.
Since $\bk$ is perfect, every constant is a $p$th power, so $j(E)$ is not constant.
Supersingular $j$-invariants lie in $\FF_{p^2}\subset\bk$, so the geometric generic fiber of $\pi$ is ordinary.
By \Cref{prop:existence of seperable rational multisection}, the geometric generic Artin--Tate twist $X_a\to\PP^1_\Omega$ has a separable smooth rational multisection of degree $p$.
Moreover $E(\bk(t)^{\sep})[p]=0$ by \Cref{sweep:no-torsion}.
Hence $J$ is unirational by \Cref{sweep:criterion}.
\end{proof}

\section{Proof of main results}\label{sec:proof}
\subsection{Frame lattices in Artin invariant nine}\label{subsec:frame-lattices}
In this subsection $p$ is an odd prime with $p\not\equiv1\pmod3$.
We construct an even positive definite lattice of rank $20$ with discriminant group $(\ZZ/p)^{18}$ and with very few roots, and realize it, up to sign, as the frame lattice of a Jacobian genus-one fibration on every supersingular $K3$ surface of Artin invariant nine.
These fibrations are used in \Cref{pencils:large} for $p\equiv2\pmod3$ and in \Cref{pencils:three} for $p=3$.
In the next two lemmas lattices are positive definite, a root is a vector of norm two, and $M(n)$ denotes a lattice $M$ with its form multiplied by $n$.

The lattice is obtained by Kneser's method of neighbors \cite{Kneser}; see \cite{Voight} for a survey.
Let $\Lambda_0$ be an even lattice, and let $\ell$ be a prime not dividing $\det\Lambda_0$.
An even lattice $\Lambda\subset\Lambda_0\otimes\QQ$ is an \emph{$\ell$-neighbor} of $\Lambda_0$ if $\Lambda\cap\Lambda_0$ has index $\ell$ in both $\Lambda$ and $\Lambda_0$ \cite[Definition~3.10]{Voight}.
The $\ell$-neighbors of $\Lambda_0$ are the lattices
\[
 \Lambda_0(\ell,z)=\{x\in\Lambda_0 \mid (x,z)\in\ell\ZZ\}+\ZZ\,z/\ell,\qquad z\in\Lambda_0\smallsetminus\ell\Lambda_0,\quad 2\ell^2\mid z^2,
\]
and they lie in the genus of $\Lambda_0$ \cite[Proposition~3.13 and Lemma~3.12]{Voight}; see also \cite[\S2]{ScharlauHemkemeier}.
We only use the case $\ell=2$, and the following lemma records what we need about it, with a direct proof.

\begin{lemma}\label{pencils:neighbor}
Let $\Lambda_0$ be an even lattice, and let $z\in\Lambda_0$ be a vector with $8\mid z^2$ such that $(x,z)$ is odd for some $x\in\Lambda_0$.
Put
\[
 \Lambda_0'=\{x\in\Lambda_0 \mid (x,z)\in2\ZZ\},\qquad \Lambda=\Lambda_0'+\ZZ\,z/2.
\]
Then $\Lambda$ is an even lattice, $\Lambda\cap\Lambda_0=\Lambda_0'$ has index two in $\Lambda_0$ and in $\Lambda$, and $\det\Lambda=\det\Lambda_0$.
If the discriminant group of $\Lambda_0$ is $(\ZZ/p)^n$, then so is the discriminant group of $\Lambda$.
\end{lemma}
We call $\Lambda$ the \emph{$2$-neighbor} of $\Lambda_0$ defined by $z$.
If $\det\Lambda_0$ is odd, the condition on $z$ means that $z\notin2\Lambda_0$, and $\Lambda$ is the $2$-neighbor $\Lambda_0(2,z)$ in the sense above.
\begin{proof}
The sublattice $\Lambda_0'$ is the kernel of the surjection $\Lambda_0\to\ZZ/2$, $x\mapsto(x,z)\bmod2$, so it has index two in $\Lambda_0$, and it contains $z$.
If $z/2$ were in $\Lambda_0$, then $(x,z)=2(x,z/2)$ would be even for every $x\in\Lambda_0$.
Hence $z/2\notin\Lambda_0$, so $\Lambda$ is the disjoint union of $\Lambda_0'$ and $\Lambda_0'+z/2$, and $\Lambda\cap\Lambda_0=\Lambda_0'$ has index two in $\Lambda$.
For $y\in\Lambda_0'$, the pairing $(y,z/2)$ is an integer and $(y+z/2)^2=y^2+(y,z)+z^2/4$ is even, so $\Lambda$ is an even lattice.
Its determinant is $\det\Lambda_0'/4=\det\Lambda_0$.
Finally $\Lambda\otimes\ZZ_p=\Lambda_0'\otimes\ZZ_p=\Lambda_0\otimes\ZZ_p$ because $p$ is odd, so the discriminant groups of $\Lambda$ and $\Lambda_0$ have isomorphic $p$-primary parts, and they have the same order.
\end{proof}

\begin{lemma}\label{pencils:frame}
There is an even lattice $\Lambda$ of rank $20$ with discriminant group $(\ZZ/p)^{18}$ whose roots are
\begin{enumerate}[label=\textup{(\roman*)},leftmargin=2.2em]
  \item $\pm u$ for a single root $u$, if $p\equiv2\pmod3$;
  \item $\pm u_1$ and $\pm u_2$ for two orthogonal roots $u_1,u_2$, if $p=3$.
\end{enumerate}
In particular, the root system of $\Lambda$ is $A_1$ in case \textup{(i)} and $A_1^{\oplus2}$ in case \textup{(ii)}.
\end{lemma}
\begin{proof}
(i) Let $u,u'$ be simple roots of $A_2$, so that $(u,u')=-1$, and let $\lambda=(2u+u')/3$, a generator of $A_2^\vee/A_2\cong\ZZ/3$ with $\lambda^2=2/3$.
Let $w=(\lambda,\lambda)$ in the dual lattice of $A_2\perp A_2(p)$, and put
\[
 L=\bigl(A_2\perp A_2(p)\bigr)+\ZZ\,w.
\]
Then $w^2=2(p+1)/3$ is an even integer because $p\equiv2\pmod3$, so $L$ is an even lattice containing $A_2\perp A_2(p)$ with index three.
Hence $\det L=p^2$ and $L\otimes\ZZ_p=(A_2\perp A_2(p))\otimes\ZZ_p$, so $L$ has discriminant group $(\ZZ/p)^2$.
Let $\Lambda_0=L\perp E_8(p)^{\oplus2}$.
It is even of rank $20$, with discriminant group $(\ZZ/p)^{18}$.

The roots of $\Lambda_0$ are the six roots of the first summand $A_2$.
Indeed, let $x$ be a root of $\Lambda_0$.
Its components in $E_8(p)^{\oplus2}$ vanish, because $E_8(p)$ has minimum $2p>2$, so $x\in L$.
Write $x=(a+n\lambda,b+n\lambda)$ with $a,b\in A_2$ and $n\in\ZZ$.
The nonzero vectors of $A_2^\vee$ have norm at least $2/3$, so a nonzero second component would have norm at least $2p/3>2$ in $A_2(p)\otimes\QQ$.
Hence $b+n\lambda=0$, which forces $3\mid n$, and $x$ lies in the first summand $A_2$.

Let $e\in E_8$ be a vector of norm six if $p\equiv1\pmod4$ and of norm two if $p\equiv3\pmod4$, viewed in the first summand $E_8(p)$, and put $z=u+e\in\Lambda_0$.
Then $z^2$ is $2+6p$ or $2+2p$, respectively, so $8\mid z^2$, and $(u',z)=(u',u)=-1$ is odd.
Let $\Lambda$ be the $2$-neighbor of $\Lambda_0$ defined by $z$.
By \Cref{pencils:neighbor}, it is even of rank $20$ with discriminant group $(\ZZ/p)^{18}$.
It contains $u$, because $(u,z)=2$.
Let $x$ be a root of $\Lambda$.
If $x\notin\Lambda_0$, then $x\in\Lambda_0'+z/2$, so the component of $x$ in the first summand $E_8(p)$ is $(e+2y)/2$ for some $y\in E_8$.
Norms in $2E_8$ are divisible by eight, so $e\notin2E_8$ and $e+2y$ is a nonzero vector of $E_8$.
Thus this component has norm $p\,(e+2y)^2/4\geq p/2>2$, which is impossible for a root.
Hence $x\in\Lambda\cap\Lambda_0=\Lambda_0'$, so $x$ is a root of the first summand $A_2$, and $(x,u)=(x,z)$ is even.
A root of $A_2$ other than $\pm u$ has pairing $\pm1$ with $u$, so $x=\pm u$.

(ii) Let $\Lambda_0=A_2^{\oplus2}\perp E_8(3)^{\oplus2}$.
It is even of rank $20$, with discriminant group $(\ZZ/3)^{18}$, and its roots are those of the two copies of $A_2$, because $E_8(3)$ has minimum $6$.
Choose simple roots $u_1,u_1'$ of the first copy of $A_2$ and a root $u_2$ of the second.
Let $e\in E_8$ be a vector of norm four, viewed in the first summand $E_8(3)$, and put $z=u_1+u_2+e$.
Then $z^2=2+2+3\cdot4=16$, and $(u_1',z)=(u_1',u_1)=-1$ is odd.
Let $\Lambda$ be the $2$-neighbor of $\Lambda_0$ defined by $z$.
By \Cref{pencils:neighbor}, it is even of rank $20$ with discriminant group $(\ZZ/3)^{18}$.
It contains $u_1$ and $u_2$, because $(u_i,z)=2$.
Let $x$ be a root of $\Lambda$.
If $x\notin\Lambda_0$, then $x\in\Lambda_0'+z/2$, so its component in the first summand $E_8(3)$ is $(e+2y)/2$ for some $y\in E_8$.
Norms in $2E_8$ are divisible by eight, so $e\notin2E_8$ and $e+2y\neq0$.
Since $e^2=4$, the norm $(e+2y)^2=4+4(e,y)+4y^2$ in $E_8$ is divisible by four, so this component has norm $3(e+2y)^2/4\geq3$ in $E_8(3)$, which is impossible for a root.
Hence $x\in\Lambda\cap\Lambda_0=\Lambda_0'$, so $x$ is a root of one copy of $A_2$, and $(x,u_1)+(x,u_2)=(x,z)$ is even.
A root of $A_2$ other than $\pm u_i$ has pairing $\pm1$ with $u_i$, so $x\in\{\pm u_1,\pm u_2\}$.
\end{proof}

\begin{corollary}\label{pencils:frame-fibration}
Let $J$ be a supersingular $K3$ surface with $\sigma(J)=9$ in characteristic $p$, and let $\Lambda$ be the lattice of \Cref{pencils:frame}.
Then $\NS(J)\cong U\perp\Lambda(-1)$, and $J$ has a Jacobian genus-one fibration $\pi\colon J\to\PP^1_\bk$ whose frame lattice is isometric to $\Lambda(-1)$.
Its frame root system is $A_1$ if $p\equiv2\pmod3$, and $A_1^{\oplus2}$ if $p=3$.
\end{corollary}
\begin{proof}
The lattice $U\perp\Lambda(-1)$ is even of signature $(1,21)$ with discriminant group $(\ZZ/p)^{18}$, so $\NS(J)\cong U\perp\Lambda(-1)$ by \Cref{thm:RS}.
As in the proof of \Cref{prelim:jacobian}, the summand $U$ gives a Jacobian genus-one fibration $\pi\colon J\to\PP^1_\bk$ whose fiber class $f$ lies in the image of $U$ under an isometry of $\NS(J)$.
Hence its frame lattice $f^\perp/\ZZ f$ is isometric to the orthogonal complement of that image, that is, to $\Lambda(-1)$.
\end{proof}

\subsection{Ordinary Jacobian fibrations in characteristic at least five}
We now show the existence of an elliptic fibration as required in \Cref{cor:criterion p>5}.
Together with \Cref{prop:AT twistor family}\textup{(ii)}, this gives the unirationality of supersingular $K3$ surfaces with $\sigma \leq 9$ for $p \geq 5$; see \Cref{subsec:proof-main}.
\begin{proposition}\label{pencils:large}
Suppose that $p\geq5$.
Every supersingular $K3$ surface of Artin invariant nine has a Jacobian elliptic fibration whose generic elliptic curve $E/\bk(t)$ satisfies $j(E)\notin\bk(t)^p$. Moreover,
\begin{enumerate}
  \item If $p \equiv 1 \pmod 3$, then every Jacobian elliptic fibration has this property.
  \item If $p \equiv 2 \pmod 3$, then there exists such a fibration with frame root system exactly $A_1$. In particular, it has precisely one reducible fiber, of type $\KI_2$ or $\KIII$, and its Mordell--Weil group has rank $19$.
\end{enumerate}
\end{proposition}
\begin{proof}
Let $J$ be a supersingular $K3$ surface with $\sigma(J)=9$.
Suppose first that $p\equiv1\pmod3$.
By \Cref{prelim:jacobian}, $J$ has a Jacobian genus-one fibration, which is elliptic because $p\geq5$.
Let $\pi$ be any Jacobian elliptic fibration on $J$.
Assume that $j(E)$ is constant.
Then $j(E)=0$ by \Cref{cor:Kodaira-j}, so $E$ has a Weierstrass equation $y^2=x^3+a_6$ with $a_6\in\bk(t)$.
For a primitive cube root of unity $\zeta\in\bk$, the automorphism $(x,y)\mapsto(\zeta x,y)$ of $E$ over $\bk(t)$ extends to an automorphism $\iota$ of $J$ over $\PP^1_\bk$, by uniqueness of relatively minimal regular models.
A nonzero global two-form on $J$ restricts on the generic fiber to $h\,dt\wedge dx/y$ for some $h\in\bk(t)$, since $\HH^0(E,\Omega^1_{E/\bk(t)})=\bk(t)\,dx/y$.
Hence $\iota$ multiplies the global two-form by $\zeta$, and $3\mid p^9+1$ by \cite[Proposition~3.3]{Jang}.
As $p^9\equiv p\pmod3$, this contradicts $p\equiv1\pmod3$.
Thus $j(E)$ is nonconstant, and $j(E)\notin\bk(t)^p$ by \Cref{cor:Kodaira-j}.

Suppose now that $p\equiv2\pmod3$.
Let $\pi\colon J\to\PP^1_\bk$ be the Jacobian genus-one fibration of \Cref{pencils:frame-fibration}, with fiber class $f$ and frame root system $A_1$; it is elliptic because $p\geq5$.
By \Cref{lem:Kodaira}, exactly one fiber of $\pi$ is reducible, and it is of type $\KI_2$ or $\KIII$.
By \Cref{cor:Kodaira-j}, $j(E)\notin\bk(t)^p$.
The vertical classes form the group $\ZZ f\oplus R_f$ of rank two, so the Shioda--Tate formula \cite[Corollary~6.13]{SS} gives $\rank E(\bk(t))=22-2-1=19$.
\end{proof}

\subsection{Jacobian elliptic fibrations in characteristic three}\label{subsec:char-three}
In characteristic three the condition \eqref{eq:star} can fail, since fibers of type $\KI_3$ and $\KIV$ have component group $\ZZ/3$ (\Cref{lem:Kodaira}).
Moreover, the generic fiber of a Jacobian fibration may be geometrically supersingular, and then \Cref{prelim:no-inseparable} does not apply.
We therefore use the fibration of \Cref{pencils:frame-fibration}, whose frame root system $A_1^{\oplus2}$ excludes such fibers.
In the supersingular case we find two distinct smooth rational trisections, that is, multisections of degree $3$ that are smooth rational curves, which are linearly equivalent on the generic fiber; they cannot both be purely inseparable, because a supersingular elliptic curve has no nonzero $3$-torsion points.

\begin{proposition}\label{pencils:three}
Suppose that $p=3$.
Every supersingular $K3$ surface $J$ with $\sigma(J)=9$ has a Jacobian elliptic fibration $\pi\colon J\to\PP^1_\bk$ with frame root system $A_1^{\oplus2}$.
It has exactly two reducible fibers, each of type $\KI_2$ or $\KIII$, and all other fibers are irreducible, of type $\KI_0$, $\KI_1$ or $\KII$; in particular, $\pi$ satisfies \eqref{eq:star}.
Its generic elliptic curve $E/\bk(t)$ has Mordell--Weil rank $18$, and it is either ordinary with $j(E)\notin\bk(t)^3$, or geometrically supersingular with $j(E)=0$.
\end{proposition}
\begin{proof}
Let $\pi\colon J\to\PP^1_\bk$ be the Jacobian genus-one fibration of \Cref{pencils:frame-fibration}, with fiber class $f$ and frame root system $A_1^{\oplus2}$.
The argument in the proof of \Cref{lem:Kodaira}, which uses only Riemann--Roch and the nefness of $f$, shows that $R_f=\bigoplus_vR_v$.
Hence the trivial lattice of $\pi$, spanned by the zero section and the fiber components, has rank $2+\rank R_f=4$.
For a quasi-elliptic fibration with a section, the Picard number is the rank of the trivial lattice \cite[Section~12.8]{SS}.
As $\rho(J)=22$, the fibration $\pi$ is elliptic, and the Shioda--Tate formula \cite[Corollary~6.13]{SS} gives $\rank E(\bk(t))=22-4=18$.

By \Cref{lem:Kodaira}, $\pi$ has exactly two reducible fibers, each of type $\KI_2$ or $\KIII$, and all other fibers are of type $\KI_0$, $\KI_1$ or $\KII$.
All of these fibers are reduced, and their component groups have order one or two, so \eqref{eq:star} holds.

In characteristic three the only supersingular $j$-invariant is $0$.
If $j(E)\neq0$, then $E$ is ordinary, and $j(E)\notin\bk(t)^3$ by \Cref{cor:Kodaira-j}.
\end{proof}

\begin{proposition}\label{pencils:three-separable}
Suppose that $p=3$, and let $\pi\colon J\to\PP^1_\bk$ be a fibration as in \Cref{pencils:three} whose generic elliptic curve $E$ is geometrically supersingular.
Let $\Omega=\kab$.
Then the geometric generic Artin--Tate twist $X_a\to\PP^1_\Omega$ has a separable smooth rational multisection of degree $3$.
\end{proposition}
\begin{proof}
Put $X=X_a$ and $K=\Omega(t)$, and let $T$ be the generic fiber of $X\to\PP^1_\Omega$, a torsor under $E_K$.
By \Cref{prelim:generic-twist} and \Cref{prelim:artin-tate}, $X$ is a supersingular $K3$ surface of Artin invariant ten, and $X\to\PP^1_\Omega$ is \'etale locally isomorphic to $\pi_\Omega$.
Hence $X\to\PP^1_\Omega$ is an elliptic fibration with exactly two reducible fibers, each consisting of two smooth rational curves of multiplicity one, and all other fibers are irreducible.
Let $f$ be the fiber class, and for $i=1,2$ let $\Theta_i$ and $\Theta_i'$ be the components of the $i$th reducible fiber.
Thus 
\[\Theta_i+\Theta_i'=f,\qquad \Theta_i^2=-2,\qquad (\Theta_1,\Theta_2)=0,
\]
and these four curves are the vertical $(-2)$-curves of $X$.

We first show that
\begin{equation}\label{eq:three-squares}
  x^2\equiv(x,\Theta_1)^2+(x,\Theta_2)^2\pmod3\qquad\text{for every }x\in\NS(X).
\end{equation}
By \Cref{roots:lattice}, the form on $\NS(X)/3\NS(X)$ has rank two.
The Gram matrix of $\Theta_1,\Theta_2$ is congruent to the identity modulo three, so these classes span a nondegenerate plane in $\NS(X)/3\NS(X)$, whose orthogonal complement is the radical.
Hence $x\equiv(x,\Theta_1)\Theta_1+(x,\Theta_2)\Theta_2$ modulo the radical, and \eqref{eq:three-squares} follows.

Next we construct two distinct smooth rational trisections whose classes differ by a vertical class.
By \Cref{roots:reflection}\textup{(ii)}, there is a root $r$ with $(r,f)=3$.
Put $m_i=(r,\Theta_i)$.
As $r^2=-2\equiv1\pmod3$, exactly one of $m_1,m_2$ is prime to three by \eqref{eq:three-squares}; we number the reducible fibers so that it is $m_1$.
Then $c=(2m_1^2-m_2+1)/3$ is an integer, and
\[
  r'=r-m_1\Theta_1+\Theta_2+cf
\]
satisfies $(r',f)=3$, $r'^2=-4-4m_1^2+2m_2+6c=-2$ and $(r',\Theta_1)=3m_1$.
By \Cref{roots:reflection}\textup{(iii)}, reflections in vertical $(-2)$-curves take $r$ and $r'$ to the classes of smooth rational multisections $M$ and $M'$ of degree $3$.
The reflection in a vertical $(-2)$-curve $C$ changes a class $x$ by a multiple of $C$, and if $(x,f)=3$, it replaces $(x,\Theta_1)$ by $\pm(x,\Theta_1)$ modulo three: indeed $s_C(x)=x+(x,C)C$, the number $(C,\Theta_1)$ is $-2$, $2$, $0$, $0$ for $C=\Theta_1,\Theta_1',\Theta_2,\Theta_2'$, and $(x,\Theta_1')=3-(x,\Theta_1)$.
Therefore $[M]-[M']$ is vertical, while $(M,\Theta_1)\not\equiv0$ and $(M',\Theta_1)\equiv0$ modulo three.
In particular $M\neq M'$.

Let $D,D'\in T$ be the generic points of $M$ and $M'$, closed points of degree three.
A vertical divisor does not meet $T$ and $\Pic(X)=\NS(X)$, so $\Oh_T(D)\cong\Oh_T(D')$.
Suppose that $M$ and $M'$ are both purely inseparable over $\PP^1_\Omega$.
Then $D$ and $D'$ both have residue field $K'=K^{1/3}$, the unique purely inseparable extension of $K$ of degree three.
Since $K'\otimes_KK'$ is a local ring with residue field $K'$, there are points $x,x'\in T(K')$ with $D_{K'}=3[x]$ and $D'_{K'}=3[x']$ as divisors on $T_{K'}$.
Thus $3[x]-3[x']$ is a principal divisor, that is, $3(x-x')=0$ in $E(K')=\Pic^0(T_{K'})$.
As $E$ is geometrically supersingular, its $3$-torsion subgroup scheme is infinitesimal, so $E(K')[3]=0$ and $x=x'$.
Hence $D=D'$ and $M=M'$, a contradiction.
So one of $M$, $M'$ is not purely inseparable over $\PP^1_\Omega$, and it is separable because $3$ is prime.
\end{proof}

\begin{corollary}\label{cor:criterion p=3}
Suppose that $p=3$.
Every supersingular $K3$ surface of Artin invariant nine is unirational.
\end{corollary}
\begin{proof}
Let $J$ be such a surface, and let $\pi\colon J\to\PP^1_\bk$ be the Jacobian elliptic fibration of \Cref{pencils:three}, with generic elliptic curve $E/\bk(t)$.
It satisfies \eqref{eq:star}, and $E(\bk(t)^{\sep})[3]=0$ by \Cref{sweep:no-torsion}, since $E$ is geometrically supersingular or $j(E)\notin\bk(t)^3$.
The geometric generic Artin--Tate twist $X_a\to\PP^1_\Omega$, with $\Omega=\kab$, has a separable smooth rational multisection of degree $3$: by \Cref{prop:existence of seperable rational multisection}\textup{(i)} if $E$ is ordinary, and by \Cref{pencils:three-separable} otherwise.
Hence $J$ is unirational by \Cref{sweep:criterion}.
\end{proof}

\subsection{Proof of the main theorems}\label{subsec:proof-main}

\begin{proof}[Proof of \Cref{thm:main}]
Let $Y$ be a supersingular $K3$ surface over $\bk$ with $\sigma(Y)\leq9$.

Suppose first that $p$ is odd.
Every supersingular $K3$ surface of Artin invariant nine is unirational.
For $p\geq5$, it has the Jacobian elliptic fibration of \Cref{pencils:large}, so \Cref{cor:criterion p>5} applies.
For $p=3$, this is \Cref{cor:criterion p=3}.
Applying \Cref{prop:AT twistor family}\textup{(ii)} successively with $s=8,7,\dots,1$ shows that every supersingular $K3$ surface of Artin invariant at most nine is unirational, and in particular $Y$ is.

Suppose now that $p=2$.
By Shimada's double-plane theorem \cite[Proposition~1.3]{ShimadaRDP}, the surface $Y$ is birational to a surface $w^2=G(x,y)$ with a polynomial $G\in\bk[x,y]$.
Since $\bk$ is perfect, $G(x,y)$ is the square of a polynomial in $x^{1/2}$ and $y^{1/2}$.
Hence the function field $\bk(x,y,w)$ of $Y$ lies in the purely transcendental field $\bk(x^{1/2},y^{1/2})$, and $Y$ is dominated by a rational surface.
This argument also applies to Artin invariant ten.
\end{proof}

\begin{proof}[Proof of \Cref{thm:B}]
In characteristic two, the double-plane argument in the proof of \Cref{thm:main} also applies to Artin invariant ten.
Assume that $p$ is odd, and let $J\to\PP^1_\bk$ be the relative Jacobian of an elliptic fibration on $X$ with supersingular geometric generic fiber.
The Jacobian $J$ is a $K3$ surface with $\rho(J)=\rho(X)=22$ \cite[p.~547 and Lemma~(1.8)]{Artin}.
The proof of \cite[Proposition~4.3.16]{BLtw}, applied to the trivial gerbe on $X$ and numerical parameter $0$, identifies $J$ with the corresponding moduli space and shows that $X$ is \'etale locally isomorphic to $J$ over $\PP^1_\bk$.
The generic fiber of $X$ thus defines a Tate--Shafarevich class for the generic fiber of $J$, and $X$ is the Artin--Tate twist of $J$ associated to the corresponding Brauer class by \Cref{lem:ArtinTate,prelim:artin-tate}.
This twist is nontrivial, since the fibration on $X$ has no section by \Cref{prelim:jacobian}.
Hence \Cref{prelim:artin-tate} gives $\sigma(J)=\sigma(X)-1=9$.
By \Cref{lem:purely inseparable sandwich}, there is a purely inseparable dominant rational map
\[
  J^{(1/p^2)}\dashrightarrow X
\]
of degree $p^2$.
Inverse Frobenius twisting preserves the N\'eron--Severi lattice, so $J^{(1/p^2)}$ is a supersingular $K3$ surface of Artin invariant nine.
It is unirational by \Cref{thm:main}, and therefore so is $X$.
\end{proof}

\begin{remark}\label{rem:BL22-multisection}
For $p\geq5$, \cite[Proposition~6.10(i)]{BLpoints} provides a purely inseparable multisection on $X$ in the setting of \Cref{thm:B}, since $J\to\PP^1_\bk$ is isotrivial with supersingular generic fiber.
This multisection trivializes the generic torsor after a suitable Frobenius base change of $\PP^1_\bk$, so the pullbacks of $X$ and $J$ are birational.
The pullback of $J$ is dominated by an inverse Frobenius twist of $J$; hence the unirationality of $J$ implies that of $X$.
\end{remark}

For $p\geq5$, the hypothesis of \Cref{thm:B} can hold only when $p\equiv2\pmod3$.
The following example gives a unirational $K3$ surface of Artin invariant ten in characteristic $p\equiv1\pmod3$.

\begin{example}\label{ex:ordinary Fermat}
Let $k=\overline{\FF}_{13}$, and let
\[
  X=\{y^2=u^5z+uv^5-z^6\}\subset\PP_k(1,1,1,3).
\]
Then $X$ is a unirational $K3$ surface of Artin invariant $10$, and every elliptic fibration on $X$ has ordinary geometric generic fiber.
Moreover, up to $k$-isomorphism, $X$ is the unique supersingular $K3$ surface of Artin invariant $10$ with non-symplectic index $13^{10}+1$.

The partial derivatives of $u^5z+uv^5-z^6$ have no common zero in $\PP^2_k$, so the branch sextic is smooth and $X$ is a $K3$ surface.
On the chart $z=1$, the substitution
\[
  (u,v,y)=\bigl(-b^{10},\,c^{10}/b^2,\,a^{25}\bigr)
\]
defines a dominant rational map from the Fermat surface $a^{50}+b^{50}+c^{50}+1=0$ to $X$ \cite[\S4]{LSY}.
Since $13^{10}\equiv-1\pmod{50}$, this Fermat surface is unirational by \cite[Theorem~III]{SK}.
Hence $X$ is unirational and therefore supersingular.

Let $\zeta\in k$ be a primitive $25$th root of unity.
The automorphism $(u,v,y)\mapsto(\zeta^{20}u,\zeta v,y)$ multiplies the regular two-form $du\wedge dv/y$ by $\zeta^{21}$, which has order $25$.
Thus the non-symplectic index of $X$ is divisible by $25$, and $25\mid 13^{\sigma(X)}+1$ by \cite[Proposition~3.3]{Jang}.
Since $13$ has order $20$ modulo $25$ and $1\leq\sigma(X)\leq10$, we obtain $\sigma(X)=10$.
At Artin invariant $10$, the possible non-symplectic indices are $2$, $13^2+1$, and $13^{10}+1$.
Only the last is divisible by $25$, and it occurs for a unique $k$-isomorphism class \cite[Table~1]{index}.
\end{example}
\bibliographystyle{amsalpha}
\bibliography{references}
\end{document}